\documentclass[12pt]{article}
\usepackage[ansinew]{inputenc}
\usepackage[pdftex]{graphicx}
\usepackage{amssymb,amsmath,amscd}
\usepackage{dsfont}
\usepackage{verbatim}
\usepackage{amsmath} 
\usepackage{graphicx}
\usepackage{pstricks}
\usepackage{enumitem}
\usepackage{ntheorem}
\usepackage{amssymb}
\usepackage{float,caption}

\newtheorem{de}{Definition}[section] 
\newtheorem{nott}{Notation}
\newtheorem{ex}[de]{Example}
\newtheorem{theo}{Theorem}
\newtheorem{prop}[theo]{Proposition}  
\newtheorem{lem}{Lemma}[section]
\newtheorem{cor}[lem]{Corollary}

\newtheorem{rem}[lem]{Remark}

\setenumerate[1]{font=\bfseries,label=\Roman*.}
\setenumerate[2]{font=\itshape,label=\arabic*)}

\author{Yashar Memarian\\
        Laboratoire de Math\'ematiques d'Orsay,\\
        Univ Paris-Sud,\\
        Orsay, F-91405,\\
        CNRS,Orsay,France}

\title{On Gromov's Waist of the Sphere Theorem}
\date{}

\begin{document}
\maketitle

\begin{abstract}
The goal of this paper is to give a detailed and complete proof of M. Gromov's waist of the sphere theorem.
\end{abstract}

\section{Introduction}

\footnotetext{E-mail address: yashar.memarian@math.u-psud.fr}

In this paper we provide details of the proof of the following important theorem.

\begin{theo}[Gromov 2003, see \cite{grwst}]
\label{1.1}
Let $f : \mathbb{S}^n \to \mathbb{R}^k$ be a continuous map from the canonical unit $n$-sphere to a Euclidean space of dimension $k$ where $k \leq n$. There exists a point $z \in \mathbb{R}^k$ such that the $n$-spherical volume of the $\varepsilon$- tubular neighborhood of $f^{-1}(z)$, denoted by $f^{-1}(z)+ \varepsilon$ satisfies, for every $\varepsilon > 0$,
$$ vol_n (f^{-1}(z)+ \varepsilon) \geq vol_n (S^{n-k}+\varepsilon).$$
Here $\mathbb{S}^{n-k}$ is the $(n-k)$-equatorial sphere of $\mathbb{S}^n$.
\end{theo}

Clearly, the Min-Max quantity dealt with in Theorem \ref{1.1} (supremum of volumes of $\varepsilon$-neighborhoods of fibers, minimized over all continuous maps to $\mathbb{R}^k$) makes sense for arbitrary metric-measure spaces. Let us call it {\em $k$-waist}. It indicates how big the space is in codimension $k$. One can see the waist as a generalization for the concentration of the measure phenomenon (which corresponds to $k=1$). The generalization has a strong topological character which is absent from classical concentration. 

M. Gromov has defined other metric measurements of $k$-dimensional size: {\em $k$-widths} are quantities which describe the thickness (diameter) of the space in codimension $k$ and {\em $k$-volumes} of maps describe how big the $k$-codimensional Hausdorff measure of the fibers of a map can be.

The proof of Theorem \ref{1.1} contains lots of interesting ideas from algebraic topology and measure theory. The first one is a generalization of the classical Borsuk-Ulam theorem which produces partitions of the sphere into finitely many convex sets whose centers have 
the same image under the given continuous map $f$. Passing to a limit, one obtains 
partitions of the sphere into infinitely many lower dimensional convex sets. Gromov 
claims that one can arrange that, for a different notion of center, the centers of the 
pieces of the limiting partition have the same image under $f$. We have been unable to 
prove this along the lines indicated by Gromov (section 5.9 in \cite{grwst}). Instead, 
in section 5, we prove weaker statements which suffice to complete the proof of Theorem 1.

We found some holes in the prove of Gromov in \cite{grwst} (section $5.5$ to $5.9$) where we fix in this paper.

The scheme of proof of the main theorem \ref{1.1} follows Gromov's paper \cite{grwst} where the author tries to complete the details missing and cover the holes. 

\section{A generalisation of the Borsuk-Ulam theorem}

Let $k=n$ and $\varepsilon = \frac{\pi}{2}$ in Theorem \ref{1.1}. In other words, let $f: \mathbb{S}^n \to \mathbb{R}^n$ be a continuous map. Theorem \ref{1.1} states the existence of a $z \in \mathbb{R}^n$ such that $vol_n (f^{-1}(z)+\pi/2) \geq vol_n (\{x,-x\}+\pi/2)$. But the right hand side of the inequality is equal the total volume of the sphere, so there is no choice for $f^{-1}(z)$ but to pass through two diametrally opposite points. We see that this particular case of the waist of the sphere theorem coincides with the classical Borsuk-Ulam theorem. So it is not a big surprise that the proof of the waist theorem relies on some algebraic topology arguments {\em à la} Borsuk-Ulam. We state first the classical Borsuk-Ulam theorem and then the generalization needed for the proof of Theorem \ref{1.1}.  

\subsection{The classical Borsuk-Ulam theorem}

\begin{theo}
\label{2}
Let $f: \mathbb{S}^{n}\rightarrow \mathbb{R}^k$ ($k\leq n$) be a continuous map from the $n$-sphere to Euclidean space of dimension $k$. There exists a partition of the sphere into two hemi-spheres and a point $z \in \mathbb{R}^k$ such that $f^{-1}(z)$ passes through the centers of both hemi-spheres.
\end{theo}

\emph{Remark} : 
\\ It is clear that the centers of the two hemi-spheres are two diametrally opposite points of the sphere. We gave a slightly different formulation of the clssical Bosuk-Ulam theorem which is better adapted to the generalization we will give later on. 

\bigskip

\emph{Proof of Theorem \ref{2}}

The map $x\mapsto g(x)=f(x)-f(-x)$ is a continuous map from $\mathbb{S}^n$ to $\mathbb{R}^k$. For every $i \in \{1,\cdots,k\}$ ,
$$ g_i(x)=f_i (x)-f_i (-x) : \mathbb{S}^n \to \mathbb{R}$$
is a continuous function from $\mathbb{S}^n$ to $\mathbb{R}$. And by the definition of the map $g$ we see that 
$\forall i \in \{1,\cdots,k\}$ we have $g_i (x)=-g_i (-x)$.

The canonical action of the group $\mathbb{Z}_2$ on the sphere $\mathbb{S}^n$ consists of sending every point to his diametrally opposite point. The quotient space is real projective space $\mathbb{R}P^n$. We can define an action of the group $\mathbb{Z}_2$ on $\mathbb{R}$ such that every point $x \in \mathbb{R}$ is sent to $-x$ by the non-trivial element of $\mathbb{Z}_2$. Hence, for every $i$, the function $g_i$ is equivariant for the action of the group $\mathbb{Z}_2$. Such a function defines a continuous cross section of the tautological vector bundle over $\mathbb{R}P^n$. And so $g$ defines a continuous cross section of Whitney sum of $k$ copies of the tautological vector bundle $\gamma_n$ over $\mathbb{R}P^n$.
$$ g : \mathbb{R}P^n \rightarrow  E=\underbrace{ \gamma_n \oplus \cdots \oplus \gamma_n}_{k}.$$

What remains to prove now is the existence of a zero for the continuous cross section $g$. For this, we refer to the theory of characteristic classes of vector bundles. In our case, as we are working with the actions of the group $\mathbb{Z}_2$, it is natural to use Stiefel-Whitney classes. The following classical result will be used here and later in this paper.

\begin{lem}
\label{milno}
Let $\pi:E \rightarrow V$ be a real vector bundle of rank $k$ over a manifold $V$. If the $k$-th Stiefel-Whitney class $w_k (E) \neq 0$, then every continuous cross section $s:V \rightarrow E$ has a zero.
\end{lem}

The cohomology ring of $\mathbb{R}P^n$ with coefficients in $\mathbb{Z}_2$ is $H^{*}(\mathbb{R}P^n,\mathbb{Z}_2)= \mathbb{Z}_2 [a]/a^{n+1}$ where $a \in H^1(\mathbb{R}P^n,Z_2)$ is the generator of the first cohomology group. One of the axioms defining Stiefel-Whitney classes states  that the total Stiefel-Whitney class $w=1+w_1 +\cdots+w_n$ is multiplicative under Whitney sums,
\begin{eqnarray*}
w(\xi\oplus\eta)=w(\xi)\smile w(\eta) .
\end{eqnarray*}
An other one states that $w(\gamma_n)=1+a$, see \cite{miln}. Thus
\begin{eqnarray*}
w(E)=(1+a)^k =1+ka+\binom{k}{2}a^2+\cdots+a^k,
\end{eqnarray*} 
and $w_k (E)=a^k$. As $k\leq n$, $a^k\neq 0$. So we proved that $w_k (E) \neq 0$. Lemma \ref{milno} implies that there exists a point $x \in \mathbb{R}P^n$ such that $g(x)=0$. And the proof of the theorem follows.

\emph{Remark} :
One should think of $\mathbb{R}P^n$ as the space of unoriented partitions of the sphere into two hemi-spheres.

Other proofs of the Borsuk-Ulam theorem can be found in \cite{matou}. We gave here a proof which was the best suited to Gromov's generalization.

\subsection{The Gromov-Borsuk-Ulam theorem}

We saw in the last section that the classical Borsuk-Ulam theorem proves the existence of a fiber passing through the center of two hemi-spheres. Gromov's generalization of Borsuk-Ulam consists of constructing a partition of the sphere into geodesically convex subsets of the sphere in order that there exists a fiber passing through the center points of all the convex sets of the partition. A hemi-spheres has a natural center point. For more general convex sets, several notions of center can be used. The Gromov-Borsuk-Ulam theorem applies to a large class of notions of center. 

\begin{de}
Say a subset $S$ of the sphere $\mathbb{S}^n$ is {\em convex} if $S$ is contained in a hemi-sphere and the cone on $S$ with vertex at the origin is convex in $\mathbb{R}^{n+1}$. Let $\mathbb{O}$ be the space of all open convex subsets of $\mathbb{S}^n$. The topology on the space $\mathbb{O}$ is defined by the Hausdorff distance between convex sets. A {\em center}map is a continuous map from $\mathbb{O}$ to $\mathbb{S}^n$. 
\end{de}

\emph{Remark} The center of a convex set is not necessarily contained in the convex set itself.

From now on, until further mention, we will fix a center map $c_{.}$.

\begin{theo}[Gromov 2003]
\label{3}
Let $f: \mathbb{S}^{n} \rightarrow \mathbb{R}^k$ ($k\leq n$) be a continuous map from the $n$-sphere to Euclidean space of dimension $k$. For every $i \in  \mathbb{N}$, there exists a partition of the sphere $\mathbb{S}^n$ into $2^i$ open convex sets $\{S_i\}$ of equal volumes ($=Vol(S^n)/2^i$) and such that all the center points $c_{.}(S_i)$ of the elements of partition have the same image in $\mathbb{R}^k$.
\end{theo} 

\emph{Remark} For $i=1$ and for a convenient choice of the center map $c_{.}$, we find Theorem \ref{2}. So this theorem can be seen as a generalisation of the classical Borsuk-Ulam theorem. But even for $i=1$  this theorem tells more than the classical Borsuk-Ulam theorem as there exists an infinite choice for the center map which won't coincide with the geometrical center of hemi-spheres. 

We saw in the last section that the space of unoriented partitions of the sphere into two hemi-spheres is identified with the real projective space. But what can we say for the space of partitions of the sphere for $i \geq 2$?

 The space of partitions into $2^i$ open convex sets of the sphere is an infinite dimensional space, we will define a finite dimensional subspace of the general space of partitions which will have very satisfying topological properties and will be easy to study. This finite dimensional subspace will be sufficient for the proof of the theorem \ref{3}.

\subsection{The partition space of $S^n$}

In this section, we define in an algorithmic way, a finite dimensional space which will be a subspace of the space of partitions of the sphere into $2^i$ open convex sets, for every natural number $i$.

We consider the following algorithm :
\begin{itemize}
\item First step. Divide $S^n$ by an oriented hyperplane into two equal hemi-spheres. The halving procedure is done by choosing a unit vector $v$ in $\mathbb{R}^{n+1}$, the two hemi-spheres are $H_{v}^{+} = \{x \in \mathbb{S}^n \quad (x.v) \geq 0 \}$ and $H_{v}^{-} = \{x\in \mathbb{S}^n \quad (x.v) \leq 0 \}$. The hemi-spheres are ordered and oriented by the vector $v$.
\item Inductive step. Divide every convex set obtained in the $(i-1)$-th step of the algorithm into two convex sets by an oriented hyperplane.
\end{itemize} 
After $i$ repetitions of the above algorithm, the sphere will be partitioned into $2^i$ convex sets. Some might be empty. In order to have $2^i$ convex sets we need 
$$1+2+2^2+\cdots+2^{i-1}=2^i-1$$
hyperplanes (hyperspheres). 

\begin{de}
\label{defpart}
The {\em space of $i$-step oriented partitions of $S^n$} is
$$P_i=\underbrace{S^n \times \cdots \times S^n}_{2^i-1}.$$
The index set $\{1,\ldots,2^i -1\}$ is viewed as the set of internal nodes of a rooted binary tree of depth $i+1$. Pieces of the partition correspond to leaves of the tree. Indeed, following the downward path connecting the root to a leaf, one meets nodes, i.e. unit vectors $v_1,\ldots,v_i$, and edges which tell whether one must use $H_{v_j}^{+}$ or $H_{v_j}^{+}$. The piece is the intersection $\bigcap_{j=1}^{i}H_{j}^{\pm}$ (eventually empty).
\end{de}

Next we want to define the space of unoriented partitions. Since the partition is defined in terms of paths connecting the root to leaves in a rooted tree, automorphisms of the rooted tree will permute points of $P_i$ which define the same unoriented partition. In the last section we saw an example for the case $i=1$. Two diametrally opposite points of the sphere define the same unoriented partition into hemi-spheres. For $i \geq 1$ things are more complicated. We give here another example.
 
\begin{ex}
\label{2.2}
Let $i=2$, we consider the space $P_2=S^n \times S^n \times S^n$ of oriented partition of the sphere into $4$ convex sets. 
\end{ex}
Let $(x,y,z) \in P_2$. $x$ is the first hyperplane cutting the sphere into two equal hemi-spheres, and defines the first step of the algorithm. At the second step, $y$ cuts the hemi-sphere pointed by $x$ into two convex pieces and $z$ cuts the hemi-sphere pointed by $-x$ into two convex pieces, providing the $4$ convex pieces of the partition defined by $(x,y,z)$. Consider the point $w=(x,-y,z)$ of $P_2$, we want to compare the partition defined by this point with the partition defined by $u=(x,y,z)$. $x$ defines the same first cut in both partitions. $-y$ and $y$ define the same hyperplane and they both cut the hemi-sphere pointed by $x$. At last, the hyperplane defined by $z$ cuts the hemi-sphere pointed by $-x$. Hence the two partitions defined by $u$ and $w$ are considered as the same unoriented partition. With the same argument, we can easily check that the $8$ following points of $P_2$ define the same partition 
$$ \{(x,y,z),(x,-y,z),(x,y,-z),(x,-y,-z),(-x,z,y),(-x,-z,y)(-x,z,-y),(-x,-z,-y)\}.$$
We define the space $Q_2$ as the quotient of $P_2$ by the equivalence relation defined by identifying the $8$ points of the above set. $Q_2$ is hence the space of unoriented partitions into $4$ convex sets defined by the above algorithm. 

In the next subsections we will explore the space $Q_i$ for all $i$.

\subsection{The binary tree $T_i$}

\begin{figure}
\begin{center}
\includegraphics[width=2in]{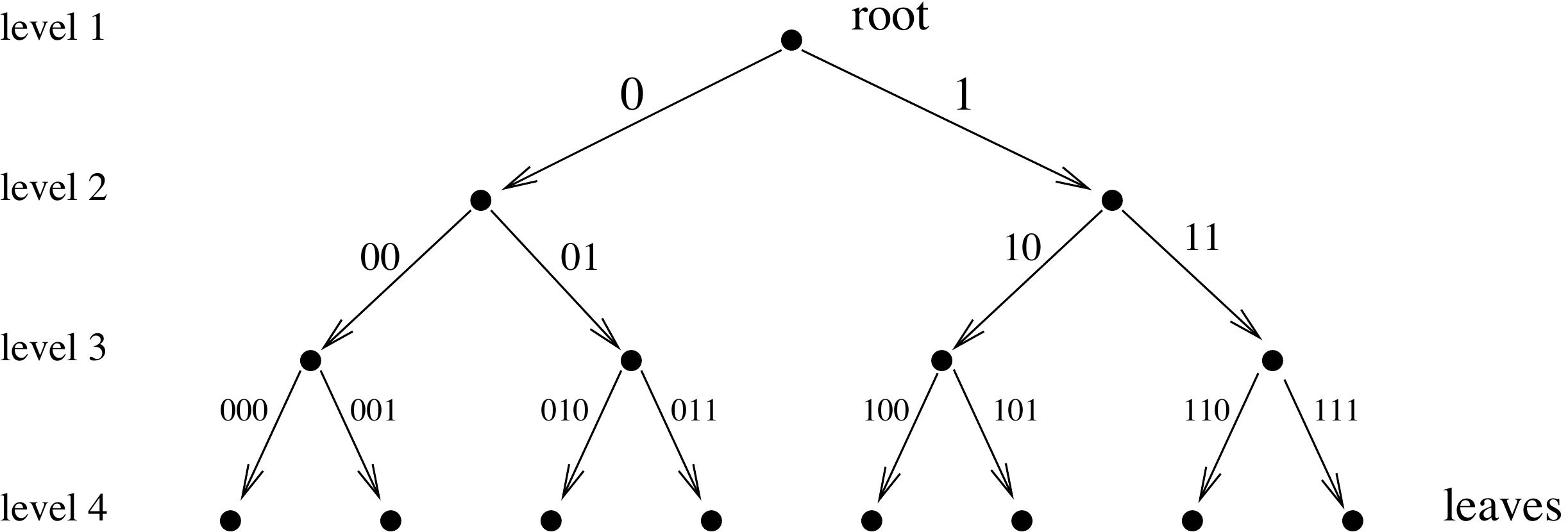}
\caption{The binary tree $T_3$}
\label{tree}
\end{center}
\end{figure}

We saw in Example \ref{2.2} that the space of oriented partitions defined as a product of some $\mathbb{S}^n$ is larger than the space of unoriented partitions. On our way to define the space of unoriented partitions, let us describe in more detail the tree structure briefly alluded to in Definition \ref{defpart}. We index the $2^i-1$ coordinates in $P_i$ by the internal nodes (i.e. vertices which are not leaves) of an oriented binary tree of depth $i$, which we denote by $T_i$. The edges are downwards oriented and indexed by strings of $0$ and $1$, as shown on Figure \ref{tree}.

Let $p=(v_n)_{n\textrm{ internal node}}\in P_i$. The unit vector $v_n$ attached to the internal node $n$ is thought of as on oriented hypersphere. To the two edges emanating from $n$ correspond hemi-spheres : the hemi-sphere to which $v_n$ points for the left edge (whose index ends with $0$), the hemi-sphere to which $-v_n$ points for the right edge (whose index ends with $1$).

\subsection{$Aut(T_i)$}

For understanding the structure of the group of automorphism $Aut(T_i)$ of the binary tree $T_i$ we need the following definition.

\begin{de}[Wreath product]
Let $G$ be a group which acts on a set $I$. Let $H$ be any group. Denote by $H^I$ the group of maps $I\to H$. The wreath product of $G$ and $H$, denoted by $G \wr H$, is the semi-direct product of the group $H^I$ by $G$,
$$G \wr H= H^I \rtimes_{\phi} G,$$
where the action $\phi$ of $G$ on $H^I$ is the left action by permuting factors,
$$ (g\cdot h)(f) = h(g^{-1}\cdot f)).$$
\end{de}

The automorphism group of a graph $\mathcal{G}$ is the set of bijections of the set of vertices such that the adjacency relationship between the vertices is respected. In other words, an automorphism of the graph $\mathcal{G}$ is a bijection $\sigma$ such that for every edge $e=uv$ where $u$ and $v$ are vertices of the graph, $\sigma(u)\sigma(v)$ is an edge of $\mathcal{G}$ (denoted by $\sigma(e)$).

\begin{lem}
for every $i \in \mathbb{N}$ we have
$$Aut(T_i)=Aut(T_{i-1}) \wr \mathbb{Z}_2 .$$
\end{lem}

\emph{Proof of the Lemma}

$G=Aut(T_{i-1})$ identifies with the subgroup of $Aut(T_i)$ which does not change the last bit in the string associated to an edge. This gives a permutation action of $Aut(T_{i-1})$ on the set $I$ of $i-1$-st level vertices of $T_i$. Note that $I$ has $2^{i-1}$ elements. It is this action which defines the wreath product. One can also view $K=(\mathbb{Z}_2)^I$ as the set of elements of $Aut(T_i)$ which fix all internal nodes. It is a normal subgroup. Indeed, any automorphism of a rooted tree permutes internal nodes. Given a leaf $\ell$ attached to an internal node $n$, denote by $b(\ell)$ denote the last bit in the string associated to the edge $n\ell$. Then $k\in(\mathbb{Z}_2)^I$ acts on leaf $\ell$ as follows : if $k(n)=0$, $k(\ell)=\ell$. Otherwise, $k(\ell)$ is the other leaf attached to $n$. In other words, $b(k(\ell))=b(\ell)+h(n)$.
 
Let $g\in G$ and $k\in K$. Then $b(g^{-1}(\ell))=b(\ell)$, $b(kg^{-1}(\ell))=b(\ell)+k(g^{-1}(n))$, $b(gkg^{-1}(\ell))=b(\ell)+k(g^{-1}(n))$. This shows that $gkg^{-1}=g\cdot k$ in $H$. Therefore the map $(k,g)\to kg\in Aut(T_i)$ defines a group homomorphism $K \rtimes_{\phi} G\to Aut(T_i)$. It is one to one, since any element of $Aut(T_i)$ coincides on internal nodes with a unique $g\in G$, and the remaining switches of leaves can be achieved by postcomposing with a unique element of $K$. Thus we get an isomorphism $K \rtimes_{\phi} G\simeq Aut(T_i)$, and the proof of the Lemma follows.

\bigskip

From Lemma \ref{2.2} we see that the automorphism group of the graph $T_i$ is formed by $i$ iterated wreath products of $\mathbb{Z}_2$ (be aware that the wreath product is not associative). And that $Aut(T_i)$ has cardinality equal to $2^{2^i-1}$.

\subsection{Unoriented partitions}

In general, If $G$ acts on a set $I$ and $H$ acts on a set $F$, $G\wr H$ acts on the set $F^I$ of maps $I\to F$ as follows. If $k\in H^I$, $z\in F^I$, $g\in G$ and $v\in I$,
\begin{eqnarray*}
kg(z)(v)=k(v)\cdot Z(g^{-1}\cdot v).
\end{eqnarray*}

\begin{de}
\label{defaction}
$Aut(T_i)$ acts on $P_i$ as follows. Elements of $Aut(T_{i-1})$ permute internal nodes, and so act by permuting the factors. If $I$ denotes the set of nodes of level $i$, elements of $K=(\mathbb{Z}_2)^I$ act on factors, with the generator indexed by $v$ acting by $x \rightarrow -x$ on the corresponding sphere factor. 

Similarly, $Aut(T_i)$ acts on $(\mathbb{R}^k)^I$.
\end{de}

Note that since the $\mathbb{Z}_2$ action on the sphere is free, the former action on $P_i$ is free.

\begin{de}
\label{defqi}
We define the {\em space of $i$-step unoriented partitions} of the sphere as the quotient space
$$ Q_i=P_i / Aut(T_i).$$
\end{de}

We have enough information to give the proof of the Gromov-Borsuk-Ulam theorem.

\subsection{Proof of Theorem \ref{3}}

Let $f$ be a continuous map from $\mathbb{S}^n$ to $\mathbb{R}^k$. Let $i \in \mathbb{N}$ be fixed and let $p \in P_i$. $p$ is a sequence of $2^i -1$ points of $\mathbb{S}^n$ that define a partition of the sphere into $2^i$ open convex sets. We represent the coordinates of $p$ by the vertices of a rooted binary tree $T_i$ of depth $i$ embedded in the plane. The $2^{i-1}$ last coordinates of $p$ are the $2^{i-1}$ hyperplanes of the last step of the algorithm. To each hyperplane $p_{.i.}$ belonging to the last $2^{i-1}$ vertices of the tree, we associate the open convex set which corresponds to the left edge outgoing from the vertex $p_{.i.}$. Hence we obtain a bijection between the $2^{i-1}$ last vertices of the tree and the left edges outgoing from each vertex. We denote this correspondance by $h_{.i.} \rightarrow S_{.i.}$ and we define the two following maps.
$$ v(h_{.i.}) = vol_n (S_{.i.}),$$
$$ \varphi(h_{.i.}) = v(h_{.i.})f(c_{.}(S_{.i.})),$$
where we remind that $c_{.}$ is the continuous center map that is supposed to be fixed.

These two maps are defined only for the hyperplanes of the $i$th step. We extend these two maps to all the hyperplanes (vertices) of $T_i$ in the following way. Let $h_{.j.}$ be a hyperplane of the $j$th step of the algorithm (a vertex of level $j$ of $T_i$). Let $T_{h_{.j.}} \subseteq T_{i}$ be the rooted binary subtree of $T_i$ whose root corresponds to $h_{.j.}$ and the edges are all the edges of $T_i$ which belongs to the subtree $T_{h_{.j.}}$. We consider the hyperplanes of the last level of the subtree $T_{h_{.j.}}$ and we define the two following maps,
$$ v(h_{.j.}) = \sum_{{h_{.i.}} \in T_{h_{.j.}}} v(h_{.i.}),$$
$$ \varphi(h_{.j.}) = \sum_{{h_{.i.}} \in T_{h_{.j.}}} \varphi(h_{.i.}).$$
Here, the sum is taken over all the vertices of level $i$ of the subtree corresponding to a vertex of level $j$.

Then we define a map $F: P_i \rightarrow (\mathbb{R}^{k+1})^{2^i-1}$ which is given by
$$ F: \{h_{.j.}\} \rightarrow \{v(h_{.j.})-v(-h_{.j.}),\varphi(h_{.j.})-\varphi(-h_{.j.})\}.$$
Since the construction only depends on the tree structure, $F$ is $Aut(T_i)$ -equivariant for the actions of $Aut(T_i)$ on $P_i$ and $(\mathbb{R}^{k+1})^{2^i-1}$. $F$ defines a continuous cross section of the vector bundle 
$$(P_i \times (\mathbb{R}^{k+1})^{2^i-1} )/Aut(T_i) \rightarrow Q_i=P_i/Aut(T_i).$$ 
The point is to show that this section vanishes. In view of Lemma \ref{2.2}, the  following characteristic class computation completes the proof of the Gromov-Borsuk-Ulam theorem. 

\begin{lem}
The top Stiefel-Whitney class of $L_i =(P_i \times (\mathbb{R}^{k+1})^{2^i-1} )/Aut(T_i)$ does not vanish.
\end{lem}

\emph{Proof of the Lemma}. 

As the action of $Aut(T_i)$ on both $P_i$ and ($\mathbb{R}^{k+1})^{2^i-1}$ is defined in an inductive way, it is natural to prove this lemma by induction. 

Since $P_i$ splits as a product $P_{i-1}\times (S^n)^{2^{i-1}}$ in a $Aut(T_{i-1})$-invariant manner, one gets a map $p_i :Q_i \to Q_{i-1}$ which is a fiber bundle with fiber $(\mathbb{R}P^n)^{2^{i-1}}$. Furthermore, since $Aut(T_{i-1})$ acts trivially on the last $2^{i-1}$ factors $\mathbb{R}^{k+1}$, on each fiber, the restriction of the bundle $L_i$ is the sum of a trivial bundle and of the bundle 
$$\underbrace{ (\gamma_n)^{k+1} \oplus \cdots \oplus (\gamma_n)^{k+1} }_{2^{i-1}}$$ 
over $(\mathbb{R}P^n)^{2^{i-1}}$. This implies that there exists a vectorbundle $\alpha_n$ on $Q_i$ whose restriction to fibers are isomorphic to $\gamma_n$, such that
\begin{eqnarray*}
L_i =p_i^* L_{i-1}\oplus \underbrace{ (\alpha_n)^{k+1} \oplus \cdots \oplus (\alpha_n)^{k+1}}_{2^{i-1}}.
\end{eqnarray*}
Thus $w(L_i)=p_i^* w(L_{i-1})\smile w(\alpha_n)^{(k+1)2^{i-1}}$. In particular, the top-dimensional components multiply,
\begin{eqnarray*}
w_{top}(L_i)=p_i^* w_{top}(L_{i-1})\smile w_1 (\alpha_n)^{(k+1)2^{i-1}}.
\end{eqnarray*}
By induction on $i$, we can assume that $w_{top}(L_{i-1})\not=0$. This implies that $w_{top}(L_{i})\not=0$.

\section{Pancakes}

Using the Gromov-Borsuk-Ulam theorem, we ideally would like to construct an infinite partition of the sphere which will have some desired properties. We know that for any continuous map from the sphere to a Euclidean space of smaller dimension, and for every natural number $i$, there exists a partition of the sphere into $2^i$ open convex sets of equal volumes and a fiber passing through the center of the convex sets of the partition. Since the volumes of the pieces of the partition tend to zero, we will have in the limit an infinite partition by convex subsets of smaller dimension. The purpose of this section is to analyse the dimension of the convex subsets when $i$ tends to infinity. How small the convex subsets can be and how can we control the dimension of the convex subsets of the partition ? In \cite{gromil}, a similar problem was considered where the sphere was sent to a two-dimensional Euclidean space and where the authors proved the existence of an infinite partition of the sphere by convex subsets of dimension $1$ using Borsuk-Ulam theorem. Here we follow the same line of ideas and by using the Gromov-Borsuk-Ulam theorem we prove the existence of an infinite partition of the sphere by convex subsets of at most dimension equal to $k$.

\begin{de}
Let $S$ be an open convex subset of $\mathbb{S}^n$, $S$ is called an $(k,\varepsilon)$-pancake if there exists a convex set $S_\pi$ of dimension $k$ such that every point of $S$ is at distance at most $\varepsilon$ from $S_\pi$.
\end{de}

\emph{Remark}  $(k,\varepsilon)$-pancakes, are used to control the dimensional size of open convex sets. For big enough $\varepsilon$ we can say that all open convex sets are pancakes. The interest of the above definition is when $\varepsilon$ is very small. In this case for a convex set to be a pancake would mean to be very close to a $k$-dimensional convex set and hence it would mean that the pancake has very small widths in $n-k$ directions orthogonal to the convex of dimension $k$.

 The typical example in Euclidean space are the rectangles, where for a rectangle of dimension $n$ to be a $k$ pancake would mean that the size of $n-k$ sides of the rectangle are very small.

Here is an improvement on Theorem \ref{3}.

\begin{theo}
\label{4}
Let $f:\mathbb{S}^n \rightarrow \mathbb{R}^k$ be a continuous map. For all $\varepsilon>0$, there exists an integer $i_0$ such that for all $i \geq i_0$ there exists a finite partition of $\mathbb{S}^n$ into $2^i$ open convex subsets such that :
\begin{enumerate}
\item Every convex subset of the partition is a $(k,\varepsilon)$-pancake.
\item The centers of all convex subsets of the partition have the same image in $\mathbb{R}^k$.
\item All convex subsets of the partition have the same volume.
\end{enumerate}
\end{theo}

\emph{Proof of the Theorem}

In the proof of Theorem \ref{3}, there was no restriction on the choice of the hyperplanes cutting the sphere. The idea of the proof of this theorem is to take a parametrized choice for the sequence of hyperplanes used to cut the sphere. 

We suppose that the sphere is cut into two equal pieces and the two center points have the same image in $\mathbb{R}^k$. Let $S_{+}$ be a hemi-sphere. We suppose that $\mathbb{S}^n$ is the unit sphere of $\mathbb{R}^{n+1}$, the boundary of the unit ball. Let $L$ be a plane of dimension $n-k-1$ passing through the origin in $\mathbb{R}^{n+1}$. Obviously $L$ intersects $S_{+}$ and the intersection locus is a half $n-k-2$-sphere. Let $L^\perp$ be the orthogonal to $L$ which we identify to a $\mathbb{R}^{k+2}$. By orthogonally projecting $\mathbb{R}^{n+1}$ onto $L^{\perp}$, every unit vector in $S^{k+1}$ defines a hyperplane (of dimension $n$), which contains $L$. So we can parametrize the hyperplanes (of dimension $n$) which contain $L$ by a sphere $S^{k+1}$.

We remember that the cutting hyperplanes of Theorem \ref{3} are indexed by their orthogonal unit vector, the idea now is to use Theorem \ref{3} by choosing every unit vector orthogonal to a hypersphere in a $S^{k+1}$. As the dimension of the range is equal to $k$, we can apply Theorem \ref{3} to the $2^i-1$ cartesian product of $S^{k+1}$ for every natural number $i$. In this case for every $i$, we obtain a partition of the sphere into $2^i$ open convex subsets of same volume, and such that in every previous step $j \leq i$, the unit vectors orthogonal to hyperplanes corresponding to this step belong to one $S^{k+1}$.

\begin{lem}
\label{knet}
For all $\varepsilon>0$, there exists an integer $N \in \mathbb{N}$ and a sequence $L_1,L_2,\ldots,L_N$ of ($n-k-1$)-dimensional planes such that for every ball of radius $\varepsilon$ in $\mathbb{S}^{k+1}$, there exists at least one $L_j$ which contains a point of that ball.
\end{lem}

\emph{Remark}  

If $k=1$, this lemma is equivalent to the existence of an $\varepsilon$-net. For $k \geq 1$ the lemma defines roughly speaking an $\varepsilon$-net in dimension $k$.

\bigskip

\emph{Proof of the Lemma}

Let $Gr(n-k-1,n+1)$ denote the Grassmannian of $(n-k-1)$-planes in $\mathbb{R}^{n+1}$. Let $L\in Gr(n-k-1,n+1)$. Let $V(L)$ be the set of all $L' \in Gr(n-k-1,n+1)$ such that $L'$ cuts the ball $B(x,\varepsilon) \cap \mathbb{S}^{k+1}$. Hence $V(L)$ is a neighbourhood of $L$ in $Gr(n-k-1,n+1)$.

The collection of $V(L)$'s defines an open covering of $Gr(n-k-1,n+1)$. By compactness, there exists a finite sub-covering and so a finite family of planes $L_1 ,\ldots L_N$ such that the $V(L_j)$ cover the Grassmannian and the proof of the lemma follows.

Lemma \ref{knet} lets us control the $l$-\emph{widths} for $l \geq k$ of pieces of the partition. Let $S_{\pi}$ be a piece of partition and let a $(k+2)$-dimensional plane passing through the origin which cuts $S_{\pi}$. By lemma \ref{knet} and the choice of the $L_i$, we can conclude that there does not exist any ball of radius $\delta$ of $\mathbb{S}^{k+1}$ in $S_{\pi} \cap S^{k+1}$. Indeed, if there exists a ball of radius $\delta$ in the intersection, then there exists a plane $L_j$ which passes through a point of this ball and hence a hyperplane $H_j$ containing $L_j$ which would cut the convex by passing through the intersecting point and this is not possible because otherwise the convex would be cut in the direction of $H_j$.

We now prove that for $\zeta$ small enough, all the $S_{\pi}$ are $(\zeta,k)$-pancakes.

\begin{lem}
For all $\zeta>0$, there exists $\varepsilon>0$ such that if $C$ is a convex set such that for every sphere $\mathbb{S}^{k+1}$, $C$ does not contain any ball of radius $\frac{\varepsilon}{4}$ of $\mathbb{S}^{k+1}$,
then $C$ is a $(\zeta,k)$-pancake.
\end{lem}

\emph{Proof of the Lemma}

By contradiction. If not, there exists a $\zeta>0$, there exists a sequence of convex sets $C_m$ which do not contain any ball of dimension $k+1$ and of radius $\varepsilon_m= \frac{1}{m}$ and which are not $(\zeta,k)$-pancakes. Let $C=\lim C_{m_j}$ where $C_{m_j}$ is a subsequence of the sequence $C_m$. Then $C$ does not contain any ball of dimension $k+1$.

Indeed $C=\lim C_{m_j}$, then for all $a, b, c \in C$ there exists $a_{m_j}, b_{m_j}, c_{m_j} \in C_{m_j}$ such that the sequences $a_{m_j} \to a$, $b_{m_j} \to b$, $c_{m_j} \to c$. By convexity, the convex hull of the three points $a_{m_j}, b_{m_j}, c_{m_j}$ : $Conv(a_{m_j},b_{m_j},c_{m_j}) \subset C_{m_j}$. But there exists $d_{m_j} \in C_{m_j}$ such that 
$$B(d_{m_j},\varepsilon/16) \subset Conv(a_{m_j},b_{m_j},c_{m_j})$$ 
and so
$$ B(d, \varepsilon/16) \subset Conv(a, b, c).$$
Hence $dim(C) \leq k$. Therefore for $m$ big enough $d_H(C_{m_j},C) \leq \zeta$ and this is a contradiction. This proof by contradiction uses Blaschke's selection principle .

This completes the proof of Theorem \ref{4}.

\section{Convexely derived measures on the sphere}

\subsection{Definition}

Remember that $\mathbb{S}^n$ is the boundary of the unit ball centered at the origin of $\mathbb{R}^{n+1}$. On $\mathbb{R}^{n+1}$ the Lebesgue measure $m_{n+1}$ is defined. We can define the (normalized) Riemannian measure on $\mathbb{S}^n$ as follows. Let $H$ be a measurable subset of $\mathbb{S}^n$. We define the set $co(H)$ by:
$$co(H)=\{\bigcup tH \vert 0 \leq t \leq 1\}.$$

The set $co(H)$ is the cone centered at the origin of $\mathbb{R}^{n+1}$ over $H$. $co(H) \subseteq \mathbb{R}^{n+1}$. We set 
$$\mu_{n}(H)=\frac{m_{n+1}(co(H))}{m_{n+1}(B_{n+1}(0,1))}.$$
$\mu_{n}$ is the normalised Riemannian measure on the sphere $\mathbb{S}^n$. 
\begin{de}
A convexely derived measure on $\mathbb{S}^n$ (resp. $\mathbb{R}^n$) is a limit of a vaguely converging sequence of probability measures of the form $\mu_i=\frac{vol|S_i}{vol(S_i)}$, where $S_i$ are open convex sets. 
\end{de}

\emph{Remark}. The support of a convexely derived measure is a convex set.   

In \cite{ale} and \cite{gromil}, the authors use concavity properties of density functions of convexely derived measures on Euclidean convex sets. Here we need also some sort of concavity properties for the density of convexely derived measures defined on convex sets of the sphere. Our approach will be to use Euclidean convex geometry by taking the cones over convex sets of the sphere and reduce spherical problems to Euclidean problems.

We begin by giving the following
\begin{de}
A real function $f$ defined on an interval of length less that $2\pi$ is called $\sin$-concave, if, when transported by a unit speed paramatrization of the unit circle, it can be extended to a $1$-homogeneous and concave function on a convex cone of $\mathbb{R}^2$.
\end{de}

This definition provides a family of example of $\sin$-concave functions. Indeed one way of obtaining a $\sin$-concave function is to consider a concave and $1$-homogeneous function on $\mathbb{R}^2$ and restrict it to $\mathbb{S}^1$. 
\begin{ex}
The linear function $f(x,y)=y$ is $1$-homogeneous and concave on $\mathbb{R}^2$. By restricting this function to the unit circle we obtain the well known function $\sin(t)$. So the sine function is $\sin$-concave.
\end{ex}
\begin{de}
A nonnegative real function $f$  is called $\sin^k$-concave if the function $f^{\frac{1}{k}}$ is $\sin$-concave.
\end{de}

The next lemma provides a familly of examples of $\sin^k$-concave functions for $k$ greater than $1$. This family will be all we need in this paper.

\begin{lem}
Let $S$ be a geodesically convex set of dimension $k$ of the sphere $\mathbb{S}^n$ with $k \leq n$. Let $\mu$ be a convexely derived measure defined on $S$ (with respect to the normalized Riemannian measure on the sphere). Then $\mu$ is a probability measure having a continuous density $f$ with respect of the canonical Riemannian measure on $\mathbb{S}^k$ restricted to $S$. Furthermore the function $f$ is $\sin^{n-k}$-concave on every geodesic arc contained in $S$.
\end{lem}

\emph{Proof of the Lemma}

Let $S_i$ be a sequence of open convex subsets of $\mathbb{S}^n$ which Hausdorff converges to $S$, where $S$ is a convex subset of dimension $k$ of the sphere. For every $i$ we define the convex cone over $S_i$ and denote it (as we saw in the beginning of this section) by $co(S_i)$. Then the sequence of open convex cones $co(S_i)$ (of dimension $(n+1)$), Hausdorff converges to the convex subset $co(S)$ (of dimension $(k+1)$). Then the sequence of normalised (probability) measures $\mu'_i=\frac{m_{n+1}|co(S_i)}{m_{n+1}(co(S_i))}$ vaguely converges to a probability measure $\mu'$ on $co(S)$. The measure $\mu'$ is convexely derived from the sequence of probability measures $\mu'_i$. We know from \cite{gromil} that the measure $\mu'$ admits a density function with respect to the $(k+1)$-dimensional Lebesgue measure and $d\mu'=F dm_{k+1}$, where $F$ is a $(n-k)$-concave function.

\begin{lem}
The measure $\mu'$ is $(n+1)$-homogeneous and the function $F$ is $(n-k)$-homogenous. Which means for every $t \in [0,1]$ and every Borel set $A$, $\mu'(tA)=t^{n+1}\mu'(A)$ and for every $x \in S$, $F(tx)=t^{n-k}F(x)$.
\end{lem}

\emph{Proof of the Lemma}

The measure $\mu'$ is convexely derived from the normalized $(n+1)$-dimensional Lebesgue measure $m_{n+1}$. $m_{n+1}$ is $(n+1)$-homogeneous and so will be for $\mu'$.
\\As $d\mu'=F dm_{k+1}$ and from the fact that $\mu'$ is $(n+1)$-homogeneous and $m_{k+1}$ is $(k+1)$-homogeneous, the function $F$ turns out to be $(n-k)$-homogeneous and the proof of the Lemma follows. 

It is then clear that the convexely derived measure $\mu$ on $S$ admits a continuous density function with respect to the canonical Riemannian measure of dimension $k$. We take two points $x$ and $y$ on $S$, take the geodesic arc $\sigma$ joining $x$ and $y$. We take the cone over $\sigma$ which is a subset of dimension $2$ of $co(S)$. We take the restriction of the function $F$ on $co(\sigma)$. We claim that the restriction of a $(n-k)$-concave function which is also $(n-k)$-homogeneous on $\sigma$ (considered as a subset of $\mathbb{S}^1$) is $\sin^{n-k}$-concave. As $F$ is $(n-k)$-concave, then $F^{1/(n-k)}$ is a concave function which is also $1$-homogeneous (as $F(tx)=t^{n-k}F(x)$ then $F^{1/(n-k)}(tx)=t F^{1/(n-k)}(x)$). Then by the previous Lemma $f^{1/(n-k)}$ is $\sin$-concave and then $f$ is $\sin^{n-k}$-concave by definition. And the proof of the main Lemma follows.
 
\subsection{More properties of $\sin$-concave functions}

\begin{lem}
\label{minima}
Let $f$ be a $\sin^k$-concave function defined on a closed interval of $\mathbb{R}$, then $f$ admits only one maximum point. Morever $f$ does not have any local minima.
\end{lem}
\emph{Proof of the Lemma}  \\
We put $g=f^{1/k}$. $g$ is $\sin$-concave. There exists a $1$-homogeneous and concave function $G$ such that $G|S=g$. Suppose $g$ has two maxima denoted by $x_1$ and $x_2$. $[x_1,x_2]$ is the segment joining these two points in $\mathbb{R}^2$. By concavity property we know that $G(\frac{x_1+x_2}{2})\geq g(x_1)=g(x_2)$. The point $x'=\frac{x_1+x_2}{2}/\vert\frac{x_1+x_1}{2}\vert \in S$. As $G$ is $1$-homogeneous we have $g(x')=G(\frac{x_1+x_2}{2})/\vert\frac{x_1+x_2}{2}\vert $ and as $\vert\frac{x_1+x_2}{2}\vert \leq 1$ then we have $g(x')\geq G(\frac{x_1+x_2}{2}) \geq g(x_1)=g(x_2)$ and this is a contradiction. Hence every $\sin^k$-function admit at most one maximum point.
\\ Suppose $g$ has a local minimum at $y$. By elementary geometry we know that there exist two points $x_1,x_2 \in S$ such that $y=\frac{x_1+x_2}{2}/\vert\frac{x_1+x_2}{2}\vert $. By the same argument as above we deduce that $g(y)\geq Min\{x_1,x_2\}$ and this is a contradiction. And the proof of the lemma follows.

\begin{lem}
\label{raccord}
Let $f$ be a continuous function defined on the interval $[-a,a]$. Assume that
\begin{itemize}
  \item $f_{|[-a,0]}$ and $f_{|[0,a]}$ are concave.
  \item the left and right derivative of $f$ at $0$ satisfy
  \begin{eqnarray*}
f'(0-)\geq f'(0+).
\end{eqnarray*}
\end{itemize}
Then $f$ is concave on the full interval $[-a,a]$.
\end{lem}

\emph{Proof of the Lemma}

Up to adding a linear function one can assume that $f'(0-)\geq 0\geq f'(0+)$. Then $f$ is nondecreasing on $[-a,0]$ and nonincreasing on $[0,a]$. For $x\in[-a,0]$, let $g_x$ be an affine function such that $g_x (x)=f(x)$ and $g_x \geq f$ on $[-a,0]$. Then $g_x$ is non increasing, thus, for $t\in[0,a]$, $g_x (t)\geq g_x (0)\geq f(0)\geq f(t)$. This shows that $g_x \geq f$ on $[-a,a]$. A similar argument applies for $x\in[0,a]$, and show that $f$ is the minimum of a family of affine functions, therefore $f$ is concave on $[-a,a]$.

\begin{lem}
\label{even}
Let $f$ be a $\sin$-concave function on an interval containing $0$, which achieves its maximum at $0$. Let $g(t)=f(|t|)$. Then $g$ is $\sin$-concave.
\end{lem}

\emph{Proof of the Lemma}

View $f$ and $g$ as functions on an arc of the unit circle in the plane containing $(1,0)$. Let $F$ and $G$ denote the $1$-homogeneous extensions of $f$ and $g$ to a plane sector $C$ containing the half line $\{(x,0)\,|\,x>0\}$. Then $G(x,y)=F(x,|y|)$ on $C$. Let $t\mapsto c(t)=(x+\alpha t,\beta t)$, $t\in[-a,a]$, be a parametrization of a line segment contained in $C$. Then $h(t)=G(c(t))$ is continuous, concave on $[-a,0]$ and $[0,a]$. Assume that $\beta>0$ and $x>0$. The left and right derivatives of $h$ at $t=0$ are equal to
\begin{eqnarray*}
h'(0-)&=&\alpha f(0)+x\beta g'(0-)=\alpha f(0)+x\beta f'(0-),\\
h'(0+)&=&\alpha f(0)+x\beta g'(0+)=\alpha f(0)-x\beta f'(0-).
\end{eqnarray*}
By assumption, $f'(0-)\geq 0$, thus $h'(0-)\geq h'(0+)$. Lemma \ref{raccord} implies that $h$ is concave. This shows that $G$ is concave, and $g$ is $\sin$-concave. 

\begin{lem}
\label{compa}
Let $0<\varepsilon<\pi/2$. Let $\tau>\varepsilon$. Let $f$ be a nonnegative $\sin^k$-concave function on $[0,\tau]$, which attains its maximum at $0$. Let $h(t)=c\cos^k (t)$ where $c$ is chosen such that $f(\varepsilon)=h(\varepsilon)$. Then 
\begin{eqnarray*}
\begin{cases}
f(x) \geq h(x)     & \text{for } x\in [0,\varepsilon], \\
f(x) \leq h(x)     & \text{for } x\in [\varepsilon,\tau].
\end{cases}
\end{eqnarray*}
In particular, $\tau\leq \pi/2$.
\end{lem}

\emph{Proof of the Lemma}

Without loss of generality, we can assume that $k=1$. Define $g(t)=f(|t|)$. View $g$ and $h$ as functions on an arc $S$ of length $\min\{\pi,2\tau\}$ of the unit circle. Let $G$ and $H$ denote the $1$-homogeneous extensions of $g$ and $h$ to the plane sector $C=co(S)$. Then $H(x,y)=cx$ on $C$. According to Lemma \ref{even}, $G$ is concave on $C$, and so is $G-H$. By construction, $G-H$ vanishes both at $p=(\cos(\varepsilon),\sin(\varepsilon))$ and at $q=(\cos(\varepsilon),-\sin(\varepsilon))$. Since $G-H$ is concave, $G-H\geq 0$ on the line segment $[p,q]$, and $G-H\leq 0$ on the remainder of $C\cap D$ where $D$ denotes the line through $p$ and $q$. Since $G-H$ is 1-homogeneous, $G-H\geq 0$ on the sector delimited by the half lines $\mathbb{R}_+ q$ and $\mathbb{R}_+ p$, and $G-H\leq 0$ on the remainder of $C$. This shows that $g\geq h$ on $[-\varepsilon,\varepsilon]$ and $g\leq h$ on $[\varepsilon,\min\{\pi/2,\tau\}]$. Assume that $\tau>\pi/2$. Then $f(\pi/2)=g(\pi/2)=h(\pi/2)=0$, so that $\pi/2$ is a local minimum of $f$. This contradicts Lemma \ref{minima}. Therefore $\tau\leq \pi/2$.

\begin{lem} 
\label{imp}
Let $\tau>0$. Let $f$ be a nonzero nonnegative $\sin^k$-concave function on $[0,\tau]$, which attains its maximum at $0$. Then $\tau\leq\pi/2$ and for all $\alpha\geq 0$ and $\varepsilon\leq\pi/2$,
$$
\frac{\int_{0}^{\min\{\varepsilon,\tau\}}f(t)\sin^{\alpha}(t)\,  dt}{\int_{0}^{\tau}f(t)\sin^{\alpha}(t) \,  dt} \geq \frac{\int_{0}^{\varepsilon}\cos^k (t)\sin^{\alpha}(t)\, dt}{\int_{0}^{\pi/2}\cos^k (t)\sin^{\alpha}(t)\, dt}.
$$
\end{lem}

\emph{Proof of the Lemma}

If $\epsilon\geq\tau$, the left hand side equals $1$, which is obviously larger than the right hand side. Otherwise, set
$$
v= \frac{\int_{0}^{\varepsilon}\cos^k(t) \sin^{\alpha}(t)\,dt}{\int_{\varepsilon}^{\pi/2}\cos^k(t) \sin^{\alpha}(t)\,dt}.
$$
Choose $c>0$ such that $h(t)=c\cos^k (t)$ satisfies $f(\varepsilon)=h(\varepsilon)$. From Lemma \ref{compa}, $\tau\leq\pi/2$, $f\geq h$ on $[0,\varepsilon]$, $f\leq h$ on $[\varepsilon,\tau]$, thus
\begin{eqnarray*}
\int_{0}^{\varepsilon}f(t) \sin^{\alpha}(t)\, dt &\geq& \int_{0}^{\varepsilon}h(t) \sin^{\alpha}(t)\, dt\\
&=& c\int_{0}^{\varepsilon}\cos^k (t)\sin^{\alpha}(t)\, dt\\
&=& cv\int_{\varepsilon}^{\pi/2}\cos^k (t) \sin^{\alpha}(t)\, dt\\
&\geq& v\int_{\varepsilon}^{\tau}h(t) \sin^{\alpha}(t)\,  dt\\
&\geq& v\int_{\varepsilon}^{\tau}f(t) \sin^{\alpha}(t)\,  dt.
\end{eqnarray*}
Thus
\begin{eqnarray*}
(1+v)\int_{0}^{\varepsilon}f(t) \sin^{\alpha}(t)\, dt\geq v\int_{0}^{\tau}f(t) \sin^{\alpha}(t)\,  dt,
\end{eqnarray*}
i.e.
\begin{eqnarray*}
\frac{\int_{0}^{\varepsilon}f(t) \sin^{\alpha}(t)\,  dt}{\int_{0}^{\tau}f(t) \sin^{\alpha}(t)\,  dt}\geq \frac{v}{1+v}=\frac{\int_{0}^{\varepsilon}\cos^k (t) \sin^{\alpha}(t)\,dt}{\int_{\varepsilon}^{\pi/2}\cos^k (t)\sin^{\alpha}(t)\,dt}.
\end{eqnarray*}

The result of Lemma \ref{imp} is very important for the estimation of the waist, as we will see in the next section.

\subsection{Lower bound for the measure of balls}

\begin{nott}
Let $\mu$ be a convexely derived measure supported on a convex set of dimension $k<n$. We denote by $M_0 (\mu)$ the unique point where its density with respect to Lebesgue $k$-dilensional measure achieves its maximum.
\end{nott}

What we need is a lower bound for $\mu(B(M_0 (\mu),\varepsilon))$. This lower bound is provided in the following Lemma.

\begin{lem} \label{ne}
Let $\mu$ be a convexely derived measure supported on a convex set $S$ of dimension $k<n$. Then 
$$\mu(B(M_0 (\mu),\varepsilon)) \geq \frac{\int_{0}^{\varepsilon}\cos^{n-k}(t)\sin^{k-1}(t) \,  dt}{\int_{0}^{\pi/2}\cos^{n-k}(t)\sin^{k-1}(t) \,  dt}. $$
\end{lem}

\emph{Proof of the Lemma}

We use polar coordinates $(t,u)\mapsto\phi(t,u)=\exp_{M_0 (\mu)}(tu)$ centered at $M_0 (\mu)$ on the $k$-sphere containing $S$: $t\in[0,\pi]$, $u\in\mathbb{S}^{k-1}$. By convexity of $S$, there exists a nonnegative function $\tau$ on $\mathbb{S}^{k-1}$ such that
\begin{eqnarray*}
\phi^{-1}(S)=\{(t,u)\,|\,0\leq t\leq \tau_u\},
\end{eqnarray*}
and
\begin{eqnarray*}
\phi^{-1}(B(M_0 (\mu),\varepsilon))=\{(t,u)\,|\,0\leq t\leq \min\{\varepsilon,\tau_u\}\},
\end{eqnarray*}
The convexely derived probability measure on $S$ is $d\mu=f\,dv$, where $dv=\sin^{k-1}(t) \,dt \,du$ and $dt$ is the Lebesgue measure on $[0,\pi]$, $du$ is the $(k-1)$-dimensional canonical Riemannian measure of $\mathbb{S}^{k-1}$. 

We shall denote abusively $f\circ\phi(t,u)$ by $f(t,u)$. Hence
$$\mu(B(M_0 (\mu),\varepsilon))=\int_{0\leq t\leq \min\{\varepsilon,\tau_u\}}f(t,u)\sin^{k-1}(t) \, dt \,du.$$
Here we can apply Lemma \ref{imp}. Let
\begin{eqnarray*}
w=\frac{\int_{0}^{\varepsilon}\cos^{n-k}(t)\sin^{k-1}(t)\, dt}{\int_{0}^{\pi/2}\cos^{n-k}(t)\sin^{k-1}(t)\, dt}.
\end{eqnarray*}
We know that for every $u \in \mathbb{S}^{k-1}$, $t\mapsto f(t,u)$ is a $\sin^{n-k}$-concave function on $[0,\tau_u]$. Therefore $\tau_u \leq \pi/2$ and
$$
\int_{0}^{\min\{\varepsilon,\tau_u\}}f(t,u)\sin^{k-1}(t)\, dt \geq w\int_{0}^{\tau_u}f(t,u)\sin^{k-1}(t) \, dt
$$
Integrating over $\mathbb{S}^{k-1}$ yields
\begin{eqnarray*}
\mu(B(M_0 (\mu),\varepsilon))&\geq& w\int_{\mathbb{S}^{k-1}}\int_{0}^{\tau_u}f(t,u)\sin^{k-1}(t) \, dt\, du\\
&=&w\mu(S)=w,
\end{eqnarray*}
since $\mu$ is a probability measure.

\begin{lem} 
\label{nee}
Let $\mathbb{S}^{n-k}$ be an equatorial $(n-k)$-dimensional sphere in $\mathbb{S}^n$ then
$$
\frac{vol_n(\mathbb{S}^{n-k}+\varepsilon)}{vol_n(\mathbb{S}^n)}=\frac{\int_{0}^{\varepsilon}\cos^{n-k}(t) \sin^{k-1}(t) \, dt}{\int_{0}^{\pi/2}\cos^{n-k}(t) \sin^{k-1}(t) \, dt}.
$$
\end{lem}

\emph{Proof of the Lemma}

Let $\mathbb{S}^{n-k}$ be an equatorial sphere. Let take the distance function from $\mathbb{S}^{n-k}$, $d(x)=d(x,\mathbb{S}^{n-k}) :\mathbb{S}^n \to \mathbb{R}$. The pushforward measure is equal $\gamma(n)\cos^{n-k}(t)\sin^{k-1} dt$, and the proof of the Lemma follows.

\section{Infinite partitions}

\begin{de}[space of convexely derived measures]
Let $\mathcal{MC}^{n}$ de\-no\-te the set of probability measures on $\mathbb{S}^n$ of the form $\mu_S =vol_{|S}/vol(S)$ where $S\subset\mathbb{S}^n$ is open and convex. The space $\mathcal{MC}$ of {\em convexely derived probability measures} on $\mathbb{S}^n$ is the vague closure of $\mathcal{MC}^{n}$.
\end{de} 

It is a compact metrizable topological space.

\begin{lem}
\label{bishop}
For all open convex sets $S\subset\mathbb{S}^n$ and all $x\in S$,
\begin{eqnarray*}
\frac{vol(S\cap B(x,r))}{vol(S)}\geq \frac{vol(B(x,r))}{vol(\mathbb{S}^n)}.
\end{eqnarray*}
\end{lem}

Proof.

Apply Bishop-Gromov's inequality in Riemannian geometry. In this special case ($\mathbb{S}^n$ has constant curvature $1$), it states that the ratio
\begin{eqnarray*}
\frac{vol(S\cap B(x,r))}{vol(B(x,r))}
\end{eqnarray*}
is a nonincreasing function of $r$. It follows that
\begin{eqnarray*}
\frac{vol(S\cap B(x,r))}{vol(B(x,r))}\geq\frac{vol(S)}{vol(\mathbb{S}^n)}.
\end{eqnarray*}

This inequality extends to all convexely derived measures, thanks to the following Lemma.

\begin{lem}
\label{radon}
{\em (See \cite{hirsh}).}
Let $\mu_i$ be a sequence of positive Radon measures on a locally compact  space $X$ which vaguely converges to a positive Radon measure $\mu$. Then for every relatively compact subset $A \subset X$ such that $\mu(\partial A)=0$,
$$\lim_{i \to \infty}\mu_i (A)=\mu(A).$$
\end{lem}

\begin{cor}
\label{bishopcdm}
For all measures $\mu\in\mathcal{MC}$ and all $x\in \mathrm{support}(\mu)$,
\begin{eqnarray*}
\mu(S\cap B(x,r))\geq \frac{vol(B(x,r))}{vol(\mathbb{S}^n)}.
\end{eqnarray*}
\end{cor}

Proof.

Let $\mu=\lim \mu_{S_{j}}$. Up to extracting a subsequence, one can assume that $S_{j}$ Hausdorff converges to a compact convex set $S$. Then $\mathrm{support}(\mu)\subset S$. Indeed, if $x\notin S$, there exists $r>0$ such that $S\cap B(x,r)=\emptyset$. Let $f$ be a continuous function on $\mathbb{S}^n$, supported in $B(x,r/2)$. Then for $j$ large enough, $S_{j}\cap B(x,r/2)=\emptyset$, $\int f\,d\mu_{S_{j}}=0$, so $\int f\,d\mu=0$, showing that $x\notin\mathrm{support}(\mu)$.

If $\mu$ is a Dirac measure, then the inequality trivially holds. Otherwise, let $x\in\mathrm{support}(\mu)$. There exist $x_j\in \mathrm{support}(\mu_j)$ such that $x_j$ tend to $x$. Since $\mu$ gives no measure to boundaries of metric balls, Lemma \ref{radon} applies, and the inequality of Lemma \ref{bishop} passes to the limit.

\begin{lem}
\label{haus}
Let $Comp(\mathbb{S}^n)$ denote the space of compact subsets of $\mathbb{S}^n$ equipped with Hausdorff distance. The map $\mathrm{support}:\mathcal{MC}\to Comp(\mathbb{S}^n)$ which maps a measure to its support is continuous.
\end{lem}

Proof.

Let $\mu_j \in\mathcal{MC}$ converge to $\mu$. One can assume that $S_j =\mathrm{support}(\mu_j)$ converge to a compact set $S$. We saw in the proof of Corollary \ref{bishopcdm} that $\mathrm{support}(\mu)\subset S$.
To prove the opposite inclusion, let us define, for $r>0$ and $x\in\mathbb{S}^n$,
\begin{eqnarray*}
f_{r,x}(y)=\begin{cases}
1 & \text{ if }d(y,x)<\frac{r}{2}, \\
2- 2\frac{d(y,x)}{r} & \text{ if }\frac{r}{2}\leq d(y,x)<r, \\
0 & \text{otherwise}.
\end{cases}
\end{eqnarray*}
Let $x\in S$. Let $x_j \in S_j$ converge to $x$. According to Lemma \ref{bishopcdm}, if  $d(x_j ,x)<r/4$,
$$
\int f_{x,r}(y)\,d\mu_{j}(y)\geq const.r^{n},
$$
i.e. $\displaystyle \int f_{x,r}\,d\mu_{j}$ does not tend to $0$. It follows that $\displaystyle \int f_{x,r}\,d\mu>0$, and $x$ belongs to $\mathrm{support}(\mu)$. This shows that $\mathrm{support}$ is a continuous map on $\mathcal{MC}$.

\medskip

The support of a convexely derived probability measure is a closed convex set, it has a dimension.

\begin{nott}
$\mathcal{MC}^{k}$ denotes the set of convexely derived probability measures whose support has dimension $k$, $\mathcal{MC}^{\leq k}=\bigcup_{\ell=0}^{k}\mathcal{MC}^{k}$, $\mathcal{MC}^{+}=\mathcal{MC}\setminus \mathcal{MC}^{0}$. For $\rho>0$, $\mathcal{MC}_{\rho}$ denotes the set of convexely derived probability measures whose support has diameter $\geq\rho$.
\end{nott}

\begin{lem}
\label{maj}
As $r$ tends to $0$, $\mu(B(x,r))$ tends to $0$ uniformly on $\mathcal{MC}_{\rho}\times\mathbb{S}^n$.
\end{lem}

Proof.

Since we deal with small radii, we can make computations as if the sphere were flat, i.e. let $\mathbb{S}^n =\mathbb{R}^n$. We can assume that $\rho$ is very small as well. Let $\mu$ be a convexely derived measure supported by a $k$-dimensional convex set $S$, let $x\in\mathbb{R}^n$ and $B=S\cap B(x,r)$. Since $S$ has diameter at least $\rho$, there is a point $y$ at distance at least $\rho/2$ of $x$. Up to a translation, we can assume that $y$ is the origin of $\mathbb{R}^k$. Let $\phi$ be the density of $\mu$. Then $\phi^{1/(n-k)}$ is concave. Thus, for $x'\in B$ and $\lambda\in]0,1[$,
\begin{eqnarray*}
\phi(\lambda x)\geq \lambda^{n-k}\phi(x).
\end{eqnarray*}
Changing variables gives
\begin{eqnarray*}
\mu(\lambda B)&=&\int_{\lambda B}\phi(z)\,dz\\
&=&\lambda^{k}\int_{B}\phi(\lambda z)\,dz\\
&\geq&\lambda^{n}\int_{B}\phi(z)\,dz\\
&=&\lambda^{n}\mu(B).
\end{eqnarray*}
If $N$ is an integer such that $N\leq \rho/4r$, then one can choose $N$ values of $\lambda$ between $1/2$ and $1$ leading to disjoint subsets $\lambda B$ of $S$, and this yields
\begin{eqnarray*}
1=\mu(S)\geq N(\frac{1}{2})^{n}\mu(B),
\end{eqnarray*}
i.e.
\begin{eqnarray*}
\mu(B)\leq 2^{n}/N \simeq \mathrm{const.}\,r/\rho.
\end{eqnarray*}

\begin{lem}
\label{equi}
The function $(\mu,x,r)\mapsto \mu(B(x,r))$ is continuous on $\mathcal{MC}_{+}\times\mathbb{S}^n \times[0,\pi/2)$.
\end{lem}

Proof.

We remind the following well known
\begin{lem}[Dini]
Let $X$ be compact, $f_j:X\to \mathbb{R}$ be a increasing (resp. decreasing) sequence of continuous functions, i.e for $i \leq i'$, $f_i\leq f_{i'}$ (resp $f_i \geq f_{i'}$). If the sequence $f_i$ is pointwise convergent then it is uniformly convergent.
\end{lem}

Fix $\rho>0$. Let $X=\mathcal{MC}_{\rho}\times \mathbb{S}^n$. Let $(\mu_j,x_j)$ converge to $(\mu,x)$. By the symmetry of the sphere, we can choose a sequence $\phi_j$ such that for every $j$, $\phi_j \in Iso(\mathbb{S}^n)$ in such a way that $\phi_j$ uniformly converges to the identity and for every $j\in \mathbb{N}$ we have $\phi_j(x_j)=x$. Hence $\mu_j(B(x_j,r))=(\phi_{j*}\mu_j)(B(x,r))$. For every $f \in C^{0}(\mathbb{S}^n)$,
$$\Vert f\circ \phi_j-f\Vert_{\infty}\underset{j \to \infty}{\longrightarrow} 0,$$
thus
$$\int_{\mathbb{S}^n}(f\circ \phi_j-f)d\mu_j \underset{j \to \infty}{\longrightarrow} 0,$$
and
\begin{eqnarray*}
\lim_{j\to \infty}\int_{\mathbb{S}^n}fd\phi_{j*}\mu_j&=&\lim_{j\to \infty}\int_{\mathbb{S}^n}f\circ\phi_j d\mu_j\\
                                                     &=&\int_{\mathbb{S}^n}f d\mu,
\end{eqnarray*}
i.e. the sequence $\phi_{j*}\mu_j$ converges vaguely to $\mu$. For every $r<\frac{\pi}{2}$ and $\mu\in \mathcal{MC}_{\rho}$, $\mu(\partial B(x,r))=0$, thus Lemma \ref{radon} applies and we conclude that $\mu_j(B(x_j,r))$ tends to $\mu(B(x,r))$. This proves that for every $r\in[0,\pi/2)$, $\mu(B(x,r))$ is a continuous function of $(\mu,x)$. 

In general, for an increasing sequence of sets $A_j$, $\mu(\bigcup A_j)=\lim_j \mu(A_j)$. This shows that for fixed $(x,\mu)$, 
\begin{eqnarray*}
\lim_{r'\to r,\,r'<r}\mu(B(x,r'))=\mu(B(x,r)),\quad \lim_{r'\to r,\,r'>r}\mu(B(x,r'))=\mu(\overline{B(x,r)}).
\end{eqnarray*} 
Again, since $\mu(\partial B(x,r))=0$, $\mu(B(x,r)$ depends continuously on $r$. Dini's Lemma implies that the function $v_{r}:(\mu,x)\mapsto \mu(B(x,r))$ varies continuously with $r$ in $C^{0}(\mathcal{MC}_{\rho}\times \mathbb{S}^n)$.

If $\mu_j\to \mu$, $x_j\to x$ and $r_j\to r$,
$$\lim_{j\to \infty}\mu_j(B(x_j,r_j))=\lim_{j\to \infty}v_{r_j}(\mu_j,x_j)=v_r(\mu,x).$$
Hence the continuity of $(\mu,x,r)\to \mu(B(x,r))$ on $\mathcal{MC}_{\rho}\times \mathbb{S}^n\times [0,\pi/2)$ and the proof of the Lemma follows.

\begin{de}[limits of finite convex partitions]
\label{toppartitions}
Let $\Pi$ be a finite convex partition of $\mathbb{S}^n$. We view it as an atomic probability measure $m(\Pi)$ on $\mathcal{MC}$ as follows: for each piece $S$ of $\Pi$, let $\mu_S =vol_{|S}/vol(S)$ be the normalized volume of $S$. Then set
\begin{eqnarray*}
m(\Pi)=\sum_{\mathrm{pieces}\,S} \frac{vol(S)}{vol(\mathbb{S}^n)}\delta_{\mu_S}.
\end{eqnarray*}
We define the {\em space of} (infinite) {\em convex partitions} $\mathcal{CP}$ as the vague closure of the image of the map $m$ in the space $\mathcal{P}(\mathcal{MC})$ of probability measures on the space of convexely derived measures. The subset $\mathcal{CP}^{\leq k}$ of convex partitions of dimension $\leq k$, consists of elements of $\mathcal{CP}$ which are supported on the subset $\mathcal{MC}^{\leq k}$ of convexely derived measures with support of dimension at most $k$.
\end{de}
Note that $\mathcal{CP}$ is compact and $\mathcal{CP}^{\leq k}$ is closed in it. Measures in the support of a convex partition can be thought of as the pieces of the partition.

\begin{lem}[desintegration formula]
\label{desint}
Let $A\subset\mathbb{S}^n$ be a set such that the intersection of $\partial A$ with every $\ell$-dimensional subsphere has vanishing $\ell$-dimensional volume, for all $\ell$, $0<\ell<n$. Let $\Pi\in \mathcal{CP}$. Assume that $\Pi(\mathcal{MC}^{0})=0$. Then
\begin{eqnarray*}
\frac{vol(A)}{vol(\mathbb{S}^n)}=\int_{\mathcal{MC}}\mu(A)\,d\Pi(\mu).
\end{eqnarray*}
\end{lem}

Proof.

The identity to be proved holds for finite partitions. According to Lemma \ref{radon}, the function $\mu\mapsto\mu(A)$ is continuous on $\mathcal{MC}^{+}$. Therefore the identity still holds for vague limits of finite partitions. This completes the proof of Lemma \ref{desint}.

\subsection{Choice of a center map}

In the previous sections, we didn't make any particular assumption about the center map. In fact the only property of this map which was used was continuity. In this section we construct a family of center maps which will lead us to the proof of the waist theorem.

\begin{de}[centers of convexely derived measures]
Let $\mu \in \mathcal{MC}$, let $r>0$. Consider the function $\mathbb{S}^n \to\mathbb{R}$, $x\mapsto v_{r,\mu} (x)=\mu(B(x,r))$. Let $M_r (S)$ be the set of points where $v_{r,\mu}$ achieves its maximum on $\mathrm{support}(\mu)$. We define the center map
$$C_{r} : \mathcal{MC} \rightarrow \mathbb{S}^n$$
by $C_{r}(S)=$ the barycenter of the convex hull of $M_r (\mu)$.

If the support of $\mu$ is $\ell$-dimensional, $0<\ell<n$, we denote by $M_0(\mu)$ the unique point where the density of $\mu$ achieves its maximum.
\end{de} 

The next Lemma states a semi-continuity property of $M_r$. 

\begin{nott}
When $A_i$, $i\in\mathbb{N}$, are subsets of a topological space, we shall denote by
\begin{eqnarray*}
\lim_{i\to\infty} A_i =\bigcap_{i}\overline{\bigcup_{j\geq i} A_j }.
\end{eqnarray*}
the set of all possible limits of subsequences $x_{i(j)} \in A_{i(j)}$.
\end{nott}

\begin{lem}
\label{semi}
Let $\mu_i$ be convexely derived measures which converge to $\mu \in\mathcal{MC}^+$. Then, for all $r>0$, 
\begin{eqnarray*}
\lim_{i\to\infty} M_r (\mu_i)\subset M_r (\mu).
\end{eqnarray*}
If follows that
\begin{eqnarray*}
\lim_{i\to\infty} \mathrm{conv.\,hull}(M_r (\mu_i))\subset \mathrm{conv.\,hull}(M_r (\mu)).
\end{eqnarray*}
\end{lem}

Proof.

Let $x\in \lim_{i\to\infty} M_r (\mu_i)$, i.e. $x=\lim_{i\to\infty}x_i$ for some $x_i \in M_r (\mu_i)$. Pick $y\in\mathrm{support}(\mu)$. Pick a sequence $y_i \in\mathrm{support}(\mu_i)$ converging to $y$. According to Lemma \ref{equi},
\begin{eqnarray*}
v_{r,\mu}(x)=\lim_{i\to\infty} v_{r,\mu_i}(x_i),\quad v_{r,\mu}(y)=\lim_{i\to\infty} v_{r,\mu_i}(y_i).
\end{eqnarray*}
Since $v_{r,\mu_i}(x_i)\geq v_{r,\mu_i}(y_i)$, we get $v_{r,\mu}(x)\geq v_{r,\mu}(y)$, showing that $x\in M_r (\mu)$.

We claim that for arbitrary compact sets $A_i \in\mathbb{S}^n$, $\lim_{i\to\infty}\mathrm{conv.\,hull}(A_i)\subset\mathrm{conv.\,hull}(\lim_{i\to\infty}A_i)$. Indeed, taking cones, it suffices to check this in Euclidean space. If $x\in\lim_{i\to\infty}\mathrm{conv.\,hull}(A_i)$, $x=\lim x_i$ with $x_i \in\mathrm{conv.\,hull}(A_i)$, then there exist $n+1$ numbers $t_{i,j}\in[0,1]$ and points $a_{i,j}\in A_i$ such that $\sum_j t_{i,j}=1$, $x_i=\sum_{j}t_{i,j}a_{i,j}$. One can assume that all sequences $i\mapsto t_{i,j}$, $a_{i,j}$ converge to $t_j$, $a_j$. Then $t_j \in[0,1]$, $\sum_j t_j =1$, $a_j \in A=\lim_{i\to\infty}A_i$ and $x=\sum_j t_j a_j \in \mathrm{conv.\,hull}(A)$. This completes the proof of Lemma \ref{semi}.

\subsection{Construction of partitions adapted to a continuous map}

\begin{de}[partitions adapted to a continuous map]
\label{f}
Let $f:\mathbb{S}^n \to \mathbb{R}^k$ be a continuous map. Let $r\geq 0$. Say a convex partition $\Pi\in\mathcal{CP}$ is {\em $r$-adapted to $f$} if there exists $z\in \mathbb{R}^k$ such that $f^{-1}(z)$ intersects the convex hull of $M_r (\mu)$ for all measures $\mu$ in the support of $\Pi$. Let 
\begin{eqnarray*}
\mathcal{F}_{r}=\{\Pi\in\mathcal{CP}\,|\,\bigcap_{\mu\in\mathrm{support}(\Pi)}f(\mathrm{conv.\,hull}(M_r (\mu)))\not=\emptyset\}
\end{eqnarray*}
denote the set of partitions which are $r$-adapted to $f$.
\end{de}

\begin{cor}
\label{ferme}
For all $r>0$, $\mathcal{F}_r$ is closed in $\mathcal{CP}$. 
\end{cor}

Proof.

If $\lim_{i\to\infty}\Pi_i =\Pi$, $\mathrm{support}(\Pi)\subset\lim_{i\to\infty}\mathrm{support}(\Pi_{i})$, i.e. every piece $\mu$ of $\Pi$ is the limit of a sequence of pieces $\mu_i$ of $\Pi_i$. By assumption, there is a $z_i \in\mathbb{R}^k$ which belongs to all $f(\mathrm{conv.\,hull}(M_r (\mu)))$, $\mu\in\mathrm{support}(\Pi_{i})$. One can assume $z_i$ converges to $z$. Then $z$ belongs to all $f(\mathrm{conv.\,hull}(M_r (\mu)))$, $\mu\in\mathrm{support}(\Pi)$. Indeed, in general, if $g$ is a continuous map and $A_i$ are subsets of a compact space, $g(\lim_{i\to\infty}A_i)=\lim_{i\to\infty}g(A_i)$. So if $\mu=\lim\mu_i$, $\mu_i \in \mathrm{support}(\Pi_{i})$,
\begin{eqnarray*}
z=\lim_{i\to\infty}z_i &\in& \lim_{i\to\infty}f(\mathrm{conv.\,hull}(M_r (\mu_i)))\\
&\subset& f\left(\lim_{i\to\infty}\mathrm{conv.\,hull}(M_r (\mu_i))\right)\\
&\subset& f(\mathrm{conv.\,hull}(M_r (\mu))),
\end{eqnarray*}
thanks to Lemma \ref{semi}.

\begin{rem}
\label{existsfr}
Theorem 3 states that for every $r>0$, $\mathcal{F}_{r}$ contains uniform atomic measures with arbitrarily many pieces. Theorem 4 produces elements of $\mathcal{F}_{r}$ whose support is contained in arbitrary thin neighborhoods of the compact subset $\mathcal{MC}^{\leq k}$. With Corollary \ref{ferme}, this gives elements in $\mathcal{F}_{r}\cap\mathcal{CP}^{\leq k}$.
\end{rem}

\subsection{Convergence of $M_r (\mu)$ as $r$ tends to $0$}

\begin{lem}
\label{convargmax}
Let $\ell<n$. For every $\ell$-dimensional convexely derived measure $\mu$,
\begin{eqnarray*}
\lim_{r\to 0}d_{H}(M_r (\mu),M_0 (\mu))=0.
\end{eqnarray*}
\end{lem}

Proof.

We prove the Lemma by contradiction. Otherwise, we get a $\delta>0$ and a
sequence of radii $r_i$ tending to $0$ such that $d_{H}(M_{r_i}(\mu),M_0 (\mu))\geq\delta$. Pick a point $x_{i} \in S$ where $v_{r_i,\mu}$ achieves its maximum and such that $d(x_{i},M_0 (\mu))\geq\delta$. Up to extracting a subsequence, we can assume that $x_{i}$ converges to $x\in S$. Then $v_{r_i,\mu}(x_i)/\alpha_k r_i^k$ converges to $\phi_{\mu}(x)$. For every $y\in S$, $v_{r_i,\mu}(y)\leq v_{r_i,\mu}(x)$ and $v_{r_i,\mu}(y)/\alpha_k r_i^k$ converges to $\phi_{\mu}(y)$. Therefore $\phi_{\mu}(y)\leq\phi_{\mu}(x)$. This shows that $\{x\}=M_0 (\mu)$, contradiction. 

\bigskip

A stronger statement will be given after the following technical lemmas.

\begin{lem}
\label{boundeddensity}
Let $\mu$ be a convexely derived measure on $\mathbb{S}^n$ whose support is a $k$-dimensional convex set $S$. Write $d\mu=\phi\,dvol_k$. Then
\begin{eqnarray*}
\max_{S}\phi\leq \frac{2^{n+1}}{vol_{k}(S)}.
\end{eqnarray*}
\end{lem}

Proof.

Replace $S$ with $C=co(S)\subset\mathbb{R}^{n+1}$, and $\phi$ by its $n-k$-homogeneous extension. Then $\phi^{1/(n-k)}$ is concave. Assume $\phi$ achieves its maximum at $x\in C$. Translate $C$ so that $x=0$. On $\frac{1}{2}C$, $\phi^{1/(n-k)}\geq\frac{1}{2}\phi^{1/(n-k)}(x)$, thus
\begin{eqnarray*}
1=\mu(S)&\geq&\int_{\frac{1}{2}C}\phi\,dvol_{k+1}\\
&\geq&\frac{1}{2^{n-k}}\phi(x)vol_{k+1}(\frac{1}{2}C)\\
&=&\frac{1}{2^{n+1}}\phi(x)vol_{k+1}(C)\\
&=&\frac{1}{2^{n+1}}\phi(x)vol_{k}(S).
\end{eqnarray*}

\begin{lem}
\label{equiconc}
Let $S$, $S_i$ be full compact convex subsets of $\mathbb{R}^{n}$ such that $S_i$ Hausdorff-converges to $S$. Let $\phi:S_i \to [0,1]$ be concave functions. Then there exists a concave function $\phi:S\to[0,1]$ and a subsequence with the following properties.
\begin{itemize}
  \item On every compact subset of the interior of $S$, $\phi_i$ converges uniformly to $\phi$.
  \item For all $x\in\partial S$ and all sequences $x_i \in S_i$ converging to $x$, 
\begin{eqnarray*}
\limsup_{i\to\infty}\phi_i (x_i) \leq \phi(x).
\end{eqnarray*}
\end{itemize}
\end{lem}

Proof.

In general, bounded concave functions $f$ on compact convex sets $\Sigma$ are locally Lipschitz, 
\begin{center}
{\em for $x\in \Sigma$ with $d(x,\partial \Sigma)=r$, and all $y\in \Sigma$, $\displaystyle |f(x)-f(y)|\leq\frac{1}{r}d(x,y)$}. 
\end{center}
Indeed, let $[x',y']$ be the intersection of $\Sigma$ with the line through $x$ and $y$, with $x'$, $x$, $y'$ and $y'$ sitting along the line in this order. Let $\ell$ be the affine function on $[x',y']$ such that $\ell(x')=f(x')$ and $\ell(x)=f(x)$. Then $f(y)\leq \ell(y)$, thus $f(y)-f(x)\leq\frac{1}{d(x',x)}|f(x)-f(x')|d(x,y)\leq\frac{1}{r}d(x,y)$. Also, let $\ell'$ be the affine function on $[x',y']$ such that $\ell'(x)=f(x)$ and $\ell'(y')=f(y')$. Then $f(y)\geq \ell'(y)$, thus $f(y)-f(x)\geq-\frac{1}{d(x,y')}|f(x)-f(y')|d(x,y)\geq-\frac{1}{r}d(x,y)$. 

This shows that on every compact subset of the interior of $S$, the sequence $f_j$ is equicontinuous, so a subsequence can be found which converges uniformly on all such compact sets to a continuous function $\phi$. Of course, $\phi$ is concave and bounded, so it extends continuously to $\partial S$. Let $x\in\partial S$ and $x_i \in S_i$ converge to $x$. Pick an interior point $x_0$ of $S$ and a second interior point $x'\not=x_0$ such that $x_0$ lies on the segment $[x',x]$. Pick $x'_i$ on the line passing through $x_0$ and $x_i$ and converging to $x'$. The Lipschitz estimate for $\phi_i$ reads
\begin{eqnarray*}
\phi_i (x_i)-\phi_i (x_0)\leq \frac{d(x_0 ,x_i)}{d(x_0 ,x'_i)}|\phi_i (x'_i)-\phi_i (x_0)|.
\end{eqnarray*}
Letting $i$ tend to infinity yields
\begin{eqnarray*}
\limsup\phi_i (x_i)\leq\phi(x_0) + \frac{d(x_0 ,x)}{d(x_0 ,x')}|\phi(x')-\phi (x_0)|.
\end{eqnarray*}
Letting $x_0$ and $x'$ tend to $x$ (while keeping $x'$, $x_0$ and $x$ aligned and $ \frac{d(x_0 ,x)}{d(x_0 ,x')}$ bounded) gives $\limsup\phi_i (x_i)\leq\phi(x)$.

\begin{lem}
\label{good}
For each $k<n$, the restriction of $(\mu,r)\mapsto d_H (M_{r}(\mu),M_0 (\mu))$ to $\mathbb{R}_+ \times \mathcal{MC}^{k}$ tends to $0$ along $\{0\} \times \mathcal{MC}^{k}$, i.e. for all $\mu\in\mathcal{MC}^{k}$,
\begin{eqnarray*}
\lim_{r\to 0,\,\mu'\to\mu,\, \mu'\in \mathcal{MC}^{k}}d_H (M_{r}(\mu),M_0 (\mu))=0.
\end{eqnarray*}
\end{lem}

Proof.

Let $\mu\in \mathcal{MC}^{k}$. Let $\mu_i$ be a sequence of $k$-dimensional convexely derived measures which converges to  $\mu$ and $r_i$ be  positive numbers tending to $0$. Let $g_i \in O(n+1)$ be a rotation mapping the support of $\mu_i$ into the $k$-sphere which contains the support of $\mu$. One can assume that $g_i$ converges to identity, and then change $\mu_i$ to $(g_{i})_{*}\mu_i$, since this does not change the convergence of centers $C_{r_i}(\mu_i)$. In other words, one can assume that all $\mu_i$ have support $S_i$ in the same $k$-sphere. Of course, $S_i$ Hausdorff-converges to the support $S$ of $\mu$. Let $\phi_i$ denote the density of $\mu_i$ with respect to $k$-dimensional volume. Since $vol_{k}(S_i)$ does not tend to $0$, $\phi_i$ are uniformly bounded, by Lemma \ref{boundeddensity}. Furthermore, on any compact convex subset $K$ of the relative interior of $S$, the $\phi_i$ are equicontinuous (this follows by the cone construction from Lemma \ref{equiconc}). Therefore one can assume that $\phi_i$ converge uniformly on compact subsets of the relative interior of $S$. Since for all $r'>0$, $v_{r',\mu_i}$ converges to $v_{r',\mu}$, the limit must be equal to the density $\phi$ of $\mu$. From Lemma \ref{equiconc}, one can assert that at boundary points $x\in\partial S$, for every sequence $x_i \in S_i$ converging to $x$, $\limsup \phi_i (x_i)\leq \phi(x)$.

We repeat the argument of Lemma \ref{convargmax}. If $M_{r_i}(\mu_i)$ does not converge to $M_0 (\mu)$, some sequence $x_i \in M_{r_i}(\mu_i)$ satisfies $d(x_{i},M_0 (\mu))\geq\delta$ for some $\delta>0$. Up to extracting a subsequence, we can assume that $x_{i}$ converges to $x\in S$. If $x\notin\partial S$, then $v_{r_i,\mu}(x_i)/\alpha_k r_i^k$ converges to $\phi(x)$. If $x\in\partial S$, $\limsup v_{r_i,\mu}(x_i)/\alpha_k r_i^k \leq \phi(x)$. For every $y\in S\setminus\partial S$, $v_{r_i,\mu}(y)\leq v_{r_i,\mu}(x)$ and $v_{r_i,\mu}(y)/\alpha_k r_i^k$ converges to $\phi(y)$. Therefore $\phi(y)\leq\phi(x)$. Since $S\setminus\partial S$ is dense in $S$, this holds for all $y\in S$, thus $\phi$ achieves its maximum at $x$, i.e. $\{x\}=M_0 (\mu)$, contradiction.

\begin{cor}
\label{unifgood}
On any compact subset of $\mathcal{MC}^{k}$, the functions 
$$
\mu\mapsto d_H (M_{r}(\mu),M_0 (\mu))
$$ 
converge uniformly to $0$ as $r$ tends to $0$.
\end{cor}

\begin{prop}
\label{goodcase}
Assume $f:\mathbb{S}^n \to \mathbb{R}^k$ is a generic smooth map. Let $r_i$ tend to $0$ and let $\Pi_i \in\mathcal{CP}^{\leq k}\cap\mathcal{F}_{r_i}$ be convex partitions of dimension $\leq k$,  $r_i$-adapted to $f$. Then, for all $\varepsilon>0$,
\begin{eqnarray*}
\max_{z\in\mathbb{R}^k}\frac{vol(f^{-1}(z)+\varepsilon)}{vol(\mathbb{S}^n)}\geq \frac{vol(\mathbb{S}^{n-k}+\varepsilon)}{vol(\mathbb{S}^n)}\limsup_{i\to\infty}\Pi_{i}(\mathcal{MC}^{k}).
\end{eqnarray*}
\end{prop}

Proof.

By assumption, for each $i$, there exists $z_i \in \mathbb{R}^k$ such that for all $\mu\in\mathrm{support}(\Pi_i)$, there exists $x_{i,\mu}\in \mathrm{conv.\,hull}(M_{r_i}(\mu))$ such that $f(x_{i,\mu})=z_i$. Let $\mathcal{K}\subset \mathcal{MC}^{k}$ be a compact set. According to Corollary \ref{unifgood} and Lemma \ref{equi}, for all $\varepsilon>0$,
\begin{eqnarray*}
\delta_i :=\sup_{\mu\in\mathcal{K}}|\mu(B(x_{i,\mu},\varepsilon))-\mu(B(M_{0}(\mu),\varepsilon))|
\end{eqnarray*}
tends to $0$. Considerations in section 5 show that for every $k$-dimensional convexely derived measure $\mu$, 
$$
\mu(B(M_{0}(\mu),\varepsilon))\geq \frac{vol(\mathbb{S}^{n-k}+\varepsilon)}{vol(\mathbb{S}^n)}.
$$
For a generic smooth map $f$, the intersection of $f^{-1}(z_i)+\varepsilon$ with $k$-dimensional convex sets has vanishing $k$-dimensional volume, so the desintegration formula applies, and
\begin{eqnarray*}
\frac{vol(f^{-1}(z_i)+\varepsilon)}{vol(\mathbb{S}^n)}
&=&\int_{\mathcal{MC}}\mu(f^{-1}(z_i)+\varepsilon)\,d\Pi_i (\mu)\\
&\geq&\int_{\mathcal{K}}\mu(B(x_{i,\mu},\varepsilon))\,d\Pi_i (\mu)\\
&\geq&\Pi_i (\mathcal{K})\frac{vol(\mathbb{S}^{n-k}+\varepsilon)}{vol(\mathbb{S}^n)}-\delta_i .
\end{eqnarray*}
Taking the supremum over all compact subsets of $\mathcal{MC}^{k}$ and then a limit as $i$ tends to infinity yields the announced inequality.

\subsection{End of the proof of Gromov's theorem}

There remains to show that convex partitions in $\mathcal{CP}^{\leq k}\cap\mathcal{F}_{r}$, $r$ small, put most of their weight on $k$-dimensional pieces. This will be proven indirectly. Pieces of dimension $<k$ may exist, but they provide a lower bound on $vol(f^{-1}(z)+r)$ which is so large, that they must have small weight. We shall need a weak concavity property of $v_{\mu,r}$, which in turn relies on the corresponding Euclidean statement.

\begin{lem}
\label{gpl}
Let $S \subset \mathbb{R}^n$ be an open convex set, $\phi$ an $m$-concave function defined on $S$. Let $\mu=\phi dvol_n$. Then the map $x\mapsto\mu(B(x,r)\cap S)$ is $m+n$-concave on $S$.
\end{lem}

Proof.

We use the following estimate (Generalized Prekopa-Leindler inequality), which can be found in \cite{led}. For $\alpha\in[-\infty,+\infty]$ and $\theta\in[0,1]$, the $\alpha$-mean of two nonnegative numbers $a$ and $b$ with weight $\theta$ is
\begin{eqnarray*}
M_{\alpha}^{(\theta)}(a,b)=(\theta a^\alpha +(1-\theta)b^\alpha)^{1/\alpha}.
\end{eqnarray*}
Let $-\frac{1}{n}\leq\alpha\leq+\infty$, $\theta\in[0,1]$, $u$, $v$, $w$ nonnegative measurable functions on $\mathbb{R}^n$ such that for all $x$, $y\in\mathbb{R}^n$,
\begin{eqnarray*}
w(\theta x+(1-\theta)y)\geq M_{\alpha}^{(\theta)}(u(x),v(y)).
\end{eqnarray*}
Let $\beta=\frac{\alpha}{1+\alpha n}$. Then
\begin{eqnarray*}
\int w\geq M_{\beta}^{(\theta)}(\int u,\int v).
\end{eqnarray*}
We apply this to restrictions of $\phi$ to balls, $u=1_{B(x,r)}\phi$, $v=1_{B(y,r)}\phi$, $w=1_{B(\theta x+(1-\theta)y,r)}\phi$. By $m$-convexity of $\phi$, the assumptions of the generalized Prekopa-Leindler inequality are satisfied with $\alpha=1/m$. Then for $\beta=\frac{1}{m+n}$,
\begin{eqnarray*}
\mu(B(\theta x+(1-\theta)y),r))\geq M_{\beta}^{(\theta)}(\mu(B(x,r)),\mu(B(y,r))),
\end{eqnarray*}
which means
\begin{eqnarray*}
\mu(B(\theta x+(1-\theta)y),r))^{\frac{1}{m+n}}\geq \theta\mu(B(x,r))^{\frac{1}{m+n}}+(1-\theta)\mu(B(y,r))^{\frac{1}{m+n}}.
\end{eqnarray*}

\begin{lem}
\label{wconcave}
The functions $v_{\mu,r}$ on $\mathbb{S}^n$ are weakly concave. In other words, there exists
a constant $c=c(n)>0$ such that for every convexely derived measure $\mu$
and every sufficiently small $r>0$, if $K\subset\mathrm{support}(\mu)$, then
\begin{eqnarray*}
\min_{\mathrm{conv}(K)}v_{\mu,\frac{r}{c}}\geq c\,\min_{K}v_{\mu,r}.
\end{eqnarray*}
\end{lem}

Proof.

Since a half-sphere is projectively equivalent with Euclidean space, it suffices to prove weak concavity when $K$ consists of 2 points.

Let $\mu$ be a $k$-dimensional convexely derived measure on $\mathbb{S}^{n}$. Denote its density by $\phi$, a $\sin^{n-k}$-concave function on the support $S$ of $\mu$. Let $\Phi$ denote the $(n-k)$-homogeneous extension of $\phi$ to the cone on $S$. This is $(n-k)$-concave. Fix a point $x_0 \in\mathbb{S}^{n}$, let $\mathbb{R}^n$ denote the tangent space of $\mathbb{S}^{n}$ at $x_0$. Denote by $\phi'$ the restriction of $\Phi$ to $\mathbb{R}^n$, and $\mu'$ the measure with density $\phi'$. Lemma \ref{gpl} implies that $x'\mapsto \mu(B(x',r))$ is $(2n-k)$-concave. This implies that for every $x'$, $y'\in\mathbb{R}^n$ and $z'$ belonging to the middle third of the line segment $[x',y']$,
\begin{eqnarray*}
\mu'(B(z',r))\geq \frac{1}{3^{2n-k}}\max\{\mu'(B(x',r)),\mu'(B(y',r))\}.
\end{eqnarray*}

The radial projection from a neighborhood $V\subset \mathbb{S}^n$ of $x_0$ to $\mathbb{R}^n$ is nearly isometric and nearly maps $\phi'$ to $\phi$. Thus there exists a constant $c_1 >0$ such that if $x$, $y\in V$ and $z$ belongs to the middle third of the geodesic segment $[x,y]$,
\begin{eqnarray*}
\mu(B(z,\frac{r}{c_1}))\geq c_1 \,\max\{\mu(B(x,r)),\mu(B(y,r))\}.
\end{eqnarray*}
Covering long segments $[x,y]$ with $N$ neighborhoods like $V$ ($N$ can be bounded independantly of $n$) provides a constant $c>0$ such that for all $z\in[x,y]$ which is not too close to the endpoints,
\begin{eqnarray*}
\mu(B(z,\frac{r}{c_1^N}))\geq c_1^N \,\max\{\mu(B(x,r)),\mu(B(y,r))\}.
\end{eqnarray*}
In particular, for $c=c_1^N$,
\begin{eqnarray*}
\mu(B(z,\frac{r}{c}))\geq c \,\min\{\mu(B(x,r)),\mu(B(y,r))\}.
\end{eqnarray*}

\begin{prop}
\label{badcase}
There exists a constant $c=c(n)>0$ such that if $\Pi$ belongs to $\mathcal{F}_{r}\cap\mathcal{CP}^{\leq k}$ for some small enough $r>0$, then for all $\ell\leq k$,
\begin{eqnarray*}
\max_{z\in\mathbb{R}^k}vol(f^{-1}(z)+\frac{r}{c})\geq c\,\sum_{\ell=0}^{k}vol(\mathbb{S}^{n-\ell}+cr)\Pi(\mathcal{MC}^{\ell}).
\end{eqnarray*}
\end{prop}

Proof.

By assumption, there exists $z\in\mathbb{R}^{k}$ such that for every measure $\mu$ in the support of $\Pi$, there exists $x\in \mathrm{conv.\,hull}(M_r (\mu))$ such that $f(x)=z$. If the support of $\mu$ is $\ell$-dimensional, Lemmata \ref{semi} and \ref{ne} give
\begin{eqnarray*}
\mu(f^{-1}(z)+\frac{r}{c})&\geq& \mu(B(x,\frac{r}{c}))\\
&=&v_{\mu,\frac{r}{c}}(x)\\
&\geq& c\,\min_{M_r (\mu)}v_{\mu,r}\\
&=&c\,\max_{\mathrm{support}(\mu)}v_{\mu,r}\\
&\geq& c\,v_{\mu,r}(M_0 (\mu)) \\
&=&c\,\mu(B(M_0 (\mu),r))\\
&\geq& c\,\frac{vol(\mathbb{S}^{n-\ell}+\rho)}{vol(\mathbb{S}^n)}. 
\end{eqnarray*}
Integrating this with respect to $\Pi$ yields
\begin{eqnarray*}
\frac{vol(f^{-1}(z)+r)}{vol(\mathbb{S}^n)}&=&\int_{\mathcal{MC}}\mu(f^{-1}(z)+r)\,d\Pi(\mu)\\
&\geq& c\,\sum_{\ell=0}^{k}\frac{vol(\mathbb{S}^{n-\ell}+r)}{vol(\mathbb{S}^n)}\Pi(\mathcal{MC}^{\ell}).
\end{eqnarray*}

\textbf{Proof of Gromov's theorem}.

At last, we prove Theorem \ref{1.1}:
{\em Let $\varepsilon>0$. Let $f:\mathbb{S}^n \to\mathbb{R}^k$ be a continuous map. Then}
\begin{eqnarray*}
\max_{z\in\mathbb{R}^k}vol(f^{-1}(z)+\varepsilon)\geq vol(\mathbb{S}^{n-k}+\varepsilon).
\end{eqnarray*}

Assume first that $f$ is smooth and generic. Then there exists a constant $W$ such that for all sufficiently small $r$,
\begin{eqnarray*}
\max_{z\in\mathbb{R}^k}vol(f^{-1}(z)+r)\leq Wr^{k}.
\end{eqnarray*}
For every $r>0$, there exists a convex partition $\Pi_r \in \mathcal{CP}^{\leq k}\cap\mathcal{F}_r$ which is $r$-adapted to $f$ (Corollary \ref{existsfr}). Proposition \ref{badcase} yields
\begin{eqnarray*}
\sum_{\ell=0}^{k}vol(\mathbb{S}^{n-\ell}+r)\Pi_r (\mathcal{MC}^{\ell})\leq \frac{1}{c}\,\max_{z\in\mathbb{R}^k}vol(f^{-1}(z)+\frac{r}{c})\leq \frac{W}{c}(\frac{r}{c})^k.
\end{eqnarray*}
As $r$ tends to $0$, this implies that for all $\ell<k$, $\Pi_r (\mathcal{MC}^{\ell})$ tends to $0$, and thus $\Pi_r (\mathcal{MC}^{k})$ tends to $1$. Letting $r$ tend to $0$ in Proposition \ref{goodcase} then shows that 
\begin{eqnarray*}
\max_{z\in\mathbb{R}^k}\frac{vol(f^{-1}(z)+\varepsilon)}{vol(\mathbb{S}^n)}\geq \frac{vol(\mathbb{S}^{n-k}+\varepsilon)}{vol(\mathbb{S}^n)}.
\end{eqnarray*}

Every continuous map $f:\mathbb{S}^n \to\mathbb{R}^k$ is a uniform limit of smooth generic maps. Hausdorff semi-continuity of $X\mapsto vol(X+\varepsilon)$ then extends the result to all continuous maps. Indeed, let the continuous map $f: \mathbb{S}^n \to \mathbb{R}^k$ of the waist theorem be fixed. Let $g_j : \mathbb{S}^n \to \mathbb{R}^k$ be a sequence of $C^{\infty}$ maps such that $\delta_{j}=\Vert g_j -f \Vert_{C^0}$ tends to $0$. For every $j$, there exists a $z_j \in \mathbb{R}^k$ such that $vol(g_j^{-1}(z_j)+\varepsilon) \geq w(\varepsilon):=vol(\mathbb{S}^{n-k}+\varepsilon)$. We know that for every $j$, $g_{j}^{-1}(z_j) \subseteq f^{-1}(B(z_j,\delta_{j}))$. Then
$$vol_n (f^{-1}(B(z_j,\delta_{j}))+\varepsilon) \geq vol_n (g_{j}^{-1}(z_j)+\varepsilon) \geq w(\varepsilon).$$
Up to extracting a subsequence, we can assume that $\{z_j\}$ converges to a point $z$. There exists a decreasing sequence $\varepsilon_j \to 0$ such that for every $j$, $\vert z-z_j \vert \leq \varepsilon_j$. Then
$$f^{-1}(B(z_j,\delta_{j}))+\varepsilon \subseteq f^{-1}(B(z,\delta_{j}+\varepsilon_j))+\varepsilon,$$
thus for all $j$
$$vol_n (f^{-1}(B(z,\delta_{j}+\varepsilon_j)+\varepsilon) \geq w(\varepsilon),$$
and by Fatou Lemma
$$vol_n (\underset{j}{\bigcap}f^{-1}(B(z,\delta_{j}+\varepsilon_j))+\varepsilon) \geq w(\varepsilon).$$
If for all $j$, $x \in f^{-1}(B(z,\delta_{j}+\varepsilon_j))+\varepsilon$, then there exists $y_j$ such that $d(x,y_j) \leq \varepsilon$ and $f(y_j) \in B(z,\delta_{j}+\varepsilon_j)$. We choose a subsequence $y_k$ which converges to $y$. By construction, $d(x,y) \leq \varepsilon$, $f(y)=z$ thus $x \in f^{-1}(z)+\varepsilon$. Hence
\begin{eqnarray*}
\underset{j}{\bigcap}f^{-1}(B(z,\delta_{j}+\varepsilon_j))+\varepsilon)\subset f^{-1}(z)+\varepsilon,
\end{eqnarray*}
and
$$vol_n (f^{-1}(z)+\varepsilon) \geq w(\varepsilon).$$

\section{Acknowledgment}
I am sincerely grateful to my phd advisor, Misha Gromov, who kindly tried to answer my extremely large amount of questions and helped me to understand(?) some (tiny) part of his mathematics. I deeply thank Pierre Pansu for his constant presence and the extraordinary amount of energy that he dedicated to me. This paper is a part of my works from the three years of my PhD studentship at Universit\'e Paris-Sud.

\bibliographystyle{plain}
\bibliography{bibwstt}

\end{document}